\documentclass[a4paper]{article}
\usepackage{amsmath,amsthm,amsfonts,thmtools}
\usepackage{amssymb}
\usepackage[cmintegrals]{newtxmath}
\usepackage{graphicx}
\usepackage{listings}
\usepackage[font=footnotesize]{caption}
\usepackage[font=scriptsize]{subcaption}
\usepackage{euscript}
\usepackage{units}
\usepackage{enumerate}
\usepackage{enumitem}
\usepackage{xcolor}
\usepackage{todonotes}
\usepackage{booktabs}
\usepackage{tikz}
\usepackage{comment}
 
\DeclareMathAlphabet{\pazocal}{OMS}{zplm}{m}{n}

\usepackage[bookmarks=true,hidelinks]{hyperref}
\usepackage{bbm}
\usepackage{multirow}

\numberwithin{equation}{section}
\newtheorem{theorem}{Theorem}[section]

\newtheorem{lemma}[theorem]{Lemma}

\newtheorem{algorithm}[theorem]{Algorithm}
\newtheorem{definition}[theorem]{Definition}

\numberwithin{equation}{section}

\theoremstyle{definition}

\newtheoremstyle{myremarkstyle}{}{}{}{}{\bfseries}{.}{ }{}
\theoremstyle{myremarkstyle}
\declaretheorem[name=Remark,qed=$\blacksquare$,numberlike=theorem]{remark}

\makeatletter
\newcommand*{\intavg}{%
  \mint@l{-}{}%
}
\newcommand*{\mint@l}[2]{%
  \@ifnextchar\limits{%
    \mint@l{#1}%
  }{%
    \@ifnextchar\nolimits{%
      \mint@l{#1}%
    }{%
      \@ifnextchar\displaylimits{%
        \mint@l{#1}%
      }{%
        \mint@s{#2}{#1}%
      }%
    }%
  }%
}
\newcommand*{\mint@s}[2]{%
  \@ifnextchar_{%
    \mint@sub{#1}{#2}%
  }{%
    \@ifnextchar^{%
      \mint@sup{#1}{#2}%
    }{%
      \mint@{#1}{#2}{}{}%
    }%
  }%
}
\def\mint@sub#1#2_#3{%
  \@ifnextchar^{%
    \mint@sub@sup{#1}{#2}{#3}%
  }{%
    \mint@{#1}{#2}{#3}{}%
  }%
}
\def\mint@sup#1#2^#3{%
  \@ifnextchar_{%
    \mint@sub@sup{#1}{#2}{#3}%
  }{%
    \mint@{#1}{#2}{}{#3}%
  }%
}
\def\mint@sub@sup#1#2#3^#4{%
  \mint@{#1}{#2}{#3}{#4}%
}
\def\mint@sup@sub#1#2#3_#4{%
  \mint@{#1}{#2}{#4}{#3}%
}
\newcommand*{\mint@}[4]{%
  \mathop{}%
  \mkern-\thinmuskip
  \mathchoice{%
    \mint@@{#1}{#2}{#3}{#4}%
        \displaystyle\textstyle\scriptstyle
  }{%
    \mint@@{#1}{#2}{#3}{#4}%
        \textstyle\scriptstyle\scriptstyle
  }{%
    \mint@@{#1}{#2}{#3}{#4}%
        \scriptstyle\scriptscriptstyle\scriptscriptstyle
  }{%
    \mint@@{#1}{#2}{#3}{#4}%
        \scriptscriptstyle\scriptscriptstyle\scriptscriptstyle
  }%
  \mkern-\thinmuskip
  \int#1%
  \ifx\\#3\\\else_{#3}\fi
  \ifx\\#4\\\else^{#4}\fi  
}
\newcommand*{\mint@@}[7]{%
  \begingroup
    \sbox0{$#5\int\m@th$}%
    \sbox2{$#5\int_{}\m@th$}%
    \dimen2=\wd0 %
    \let\mint@limits=#1\relax
    \ifx\mint@limits\relax
      \sbox4{$#5\int_{\kern1sp}^{\kern1sp}\m@th$}%
      \ifdim\wd4>\wd2 %
        \let\mint@limits=\nolimits
      \else
        \let\mint@limits=\limits
      \fi
    \fi
    \ifx\mint@limits\displaylimits
      \ifx#5\displaystyle
        \let\mint@limits=\limits
      \fi
    \fi
    \ifx\mint@limits\limits
      \sbox0{$#7#3\m@th$}%
      \sbox2{$#7#4\m@th$}%
      \ifdim\wd0>\dimen2 %
        \dimen2=\wd0 %
      \fi
      \ifdim\wd2>\dimen2 %
        \dimen2=\wd2 %
      \fi
    \fi
    \rlap{%
      $#5%
        \vcenter{%
          \hbox to\dimen2{%
            \hss
            $#6{#2}\m@th$%
            \hss
          }%
        }%
      $%
    }%
  \endgroup
}

\def\XXint#1#2#3{{\setbox0=\hbox{$#1{#2#3}{\int}$ }
		\vcenter{\hbox{$#2#3$ }}\kern-.6\wd0}}

\renewcommand{\geq}{\geqslant}
\renewcommand{\leq}{\leqslant}

\renewcommand{\epsilon}{\varepsilon}
\renewcommand{\phi}{\varphi}

\newcommand{\IP}{\mathbb{P}}

\newcommand{\R}{\mathbb{R}}
\newcommand{\N}{\mathbb{N}}

\newcommand{\U}{{\bf U}}		

\newcommand{\Prob}{\EuScript{P}}

\newcommand{\IE}{\mathbb{E}}
\newcommand{\map}{\EuScript{L}}
\newcommand{\train}{\EuScript{S}}
\newcommand{\vtrain}{\mathbf{S}}
\newcommand{\test}{\EuScript{T}}
\newcommand{\val}{\EuScript{V}}
\newcommand{\reg}{\EuScript{R}}
\newcommand{\er}{\EuScript{E}}
\newcommand{\var}{\mathbb{V}}

\newcommand{\cost}{\EuScript{C}}

\newcommand{\dtl}{\EuScript{D}}

\newcommand{\E}[1]{\mathbb{E}[#1]}

\begin{document}

\date{\today}

\title{A Multi-level procedure for enhancing accuracy of machine learning algorithms.}

\author{Kjetil O. Lye \thanks{Seminar for Applied Mathematics (SAM), D-Math \newline
  ETH Z\"urich, R\"amistrasse 101, 
  Z\"urich-8092, Switzerland}, Siddhartha Mishra \thanks{Seminar for Applied Mathematics (SAM), D-Math \newline
  ETH Z\"urich, R\"amistrasse 101, 
  Z\"urich-8092, Switzerland} and Roberto Molinaro \thanks{D-Mavt,ETH Z\"urich, R\"amistrasse 101, 
  Z\"urich-8092, Switzerland .}}

\date{\today}

\maketitle
\begin{abstract}
We propose a multi-level method to increase the accuracy of machine learning algorithms for approximating observables in scientific computing, particularly those that arise in systems modeled by differential equations. The algorithm relies on judiciously combining a large number of computationally cheap training data on coarse resolutions with a few expensive training samples on fine grid resolutions. Theoretical arguments for lowering the generalization error, based on reducing the variance of the underlying maps, are provided and numerical evidence, indicating significant gains over underlying single-level machine learning algorithms, are presented. Moreover, we also apply the multi-level algorithm in the context of forward uncertainty quantification and observe a considerable speed-up over competing algorithms. 
\end{abstract}
\section{Introduction}
A fundamental goal in scientific computing is the efficient simulation of \emph{observables}, also referred to as functionals, quantities of interest or figures of merit, of systems that arise in physics and engineering. Often, the underlying system is modelled by nonlinear partial differential equations (PDEs). A prototypical example is provided by the simulation of flows past aerospace vehicles where the observables of interest are body forces such as the lift and the drag coefficients and the underlying PDEs are the compressible Euler or Navier-Stokes equations of fluid dynamics. Other interesting examples include the run-up height for a tsunami (with the shallow water equations modeling the flow) or loads (stresses) on structures, with the underlying system being modeled by the equations of elasticity or visco-elasticity. 

Computing observables involves first solving the underlying PDE by suitable numerical methods such as finite difference, finite element, finite volume or spectral methods and then evaluating the corresponding observable by another numerical method, usually a quadrature rule. 

The computation of observables can be very expensive as solving the underlying (nonlinear) PDE, especially in three space dimensions, entails a large computational cost, even on state of the art high performance computing (HPC) platforms. 

This high computational cost is particularly evident when one considers large scale problems such as uncertainty quantification, (Bayesian) inverse problems, data assimilation or optimal control/design. All these problems are of the \emph{many query} type i.e, the underlying PDE has to be solved for a very large number of instances, each corresponding to a particular realization of the input parameter space, in order to compute the \emph{input parameters to observable} map.  Querying the computationally costly PDE solver multiple times renders these problems prohibitively expensive.  

Although many methods such as reduced order models \cite{ROMbook} have been developed to provide a \emph{surrogate} for the PDE solver in computing observables, they may not be stable or accurate enough for complex non-linear PDEs such as those modeling fluid flows. Hence, there is a pressing need for the design of fast and accurate surrogates for computing observables. 

Machine learning is a very popular field within computer science in recent years. In particular, artificial neural networks i.e, layers of units (neurons) that compose affine transformations with simple (scalar) nonlinearities, are a very effective tool in a variety of contexts. \emph{Deep learning}, based on artificial neural networks with a large number of hidden layers, is extremely successful at diverse tasks, for instance in image processing, computer vision, text and speech recognition, game intelligence and more recently in protein folding \cite{AFOL}, see \cite{DL-nat} and references therein for more applications of deep learning. A key element in \emph{supervised} deep learning is the \emph{training} of tunable parameters in the underlying neural network by (approximately) minimizing suitable \emph{loss functions} on the set of training data. The resulting very high-dimensional (non-convex) optimization problem is customarily solved with variants of the stochastic gradient descent method \cite{SG}.  

Machine learning is being increasingly used in the context of scientific computing. Given that neural networks are very powerful universal function approximators \cite{Cy1,YAR1, MSDR}, it is natural to consider the space of neural networks as an ansatz space for approximating solutions of PDEs. First proposed in \cite{Lag1} on an underlying collocation approach, it has been successfully used recently in different contexts. See \cite{Kar1,Kar2, Kar3, JR1,E1,E2,E3} and references therein. This approach appears to work quite well for problems with high regularity (smoothness) of the underlying solutions and/or if the solution of the underlying PDE possesses a representation formula in terms of integrals.

Another set of methods embed deep learning modules within existing numerical methods to improve them. Examples include solving the elliptic equations in a divergence projection step in incompressible flows \cite{INC} or learning troubled cell indicators within an RKDG code for applying limiters \cite{DR1} or recasting finite difference (volume) schemes as neural networks and training the underlying parameters to improve accuracy on coarse grids \cite{DL_SM1}.

In a recent paper \cite{LMR1}, the authors developed a deep learning algorithm to approximate observables (functionals) in computational fluid dynamics. One of the challenges, discussed in \cite{LMR1}, in using deep learning algorithms in many contexts in scientific computing, stems from estimates on the so-called \emph{generalization error} (\cite{CS1}, see also 2.26 of \cite{LMR1} or \eqref{eq:ecgen} for the definition), that measures the accuracy of the trained network on unseen inputs. Although estimating generalization error sharply is a notoriously hard problem \cite{ARORA,CS1}, an upper bound on the generalization error in the context of regression of functions, for randomly selected training data is usually in the form:
\begin{equation}
    \label{eq:genE1}
    \bar{\er}_G \sim \bar{\er}_T + \frac{U}{\sqrt{N}},
\end{equation}
with $\bar{\er}_G$ being the generalization error, $\bar{\er}_T$ the training error (see \eqref{eq:ectrain} for definition) and $N$ being the number of training samples. The detailed estimate, presented in section \ref{sec:2}, \eqref{eq:genE2}, bounds $U$ in \eqref{eq:genE1} in terms of two components, one arising from the so-called \emph{validation gap} for the neural network (to check overfitting) and the other depending on the variation (measured by the standard deviation) of the underlying function (and neural network) . The bound \eqref{eq:genE1} illustrates one of the challenges of using deep learning in the context of approximating (observables of) solutions of PDEs. As long as the upper bound $U \sim {\mathcal O} (1)$, we need a large number of training samples in order to obtain reasonably small generalization errors. Since the training samples are generated by solving PDEs, generation of a large number of training samples necessitates a very high computational cost. 

In \cite{LMR1}, the authors proposed reducing the generalization error in this \emph{data poor} regime, by selecting training points based on \emph{low-discrepancy sequences}, such as Halton or Sobol sequences which are heavily used in the context of Quasi-Monte Carlo (QMC) methods \cite{CAF1}. This approach yielded considerable increase in accuracy at the same cost, over choosing random training points and allowed the authors to obtain very low prediction errors with a few ($\mathcal{O}(100))$) training samples for problems such as predicting drag and lift for flows past airfoils. 

However, one can only prove that such an approach of using low-discrepancy sequences reduces the generalization error if the underlying function is sufficiently regular (see the recent paper \cite{MR1} for the relevant estimates). In general, observables arising in computational fluid dynamics for instance, have rather low regularity (see section 2.3.2 of \cite{LMR1}). Moreover, the approach of using low-discrepancy sequences is viable only if the dimension of the underlying input parameter space is only moderately high. 

In this paper, we propose another approach for increasing the accuracy (by reducing the generalization error) of deep neural networks for computing observables. Our approach is based on the observation that the upper bound $U$ in \eqref{eq:genE1} involves the standard deviation of the underlying observable (see \eqref{eq:genE2} for an exact statement of the estimate). We will reduce the variance (standard deviation) of the observable in order to lower the generalization error. To this end, we adapt a \emph{multi-level} or multi-resolution procedure to the context of machine learning.

Multi-level methods were introduced in the context of numerical quadrature by Heinrich in \cite{Hein1} and for numerical solutions of stochastic differential equations (SDEs) by Giles in \cite{GIL1}. They have been heavily used in recent years for uncertainty quantification in PDEs, solutions of stochastic PDEs, data assimilation and Bayesian inversion. See \cite{GIL2} for a detailed survey of multi-level Monte Carlo methods and applications and \cite{MS1,MSS1} for applications to computational fluid dynamics. We would like to point out that multi-level methods are inspired by multi-grid and multi-resolution techniques, which have been used in numerical analysis over many decades. 

The basic idea of our multi-level machine learning algorithm is to approximate the observable on a sequence of nested mesh resolutions for solving the underlying PDE. We then learn the so-called \emph{details} (differences of the observable on successive mesh resolutions), instead of the observable itself. If the underlying numerical method converges to a solution of the PDE, then the standard deviation of the details is significantly smaller than the standard deviation of the underlying observable, resulting in a smaller value of the upper bound in \eqref{eq:genE1} and allowing us to learn the details with significantly fewer training samples. By carefully balancing the standard deviation of the details with the computational cost of generating the observable at each level of resolution, we aim to reduce the overall generalization error, while keeping the cost of generating the training samples small. 

The main aim of this paper is to present this novel multi-level machine learning algorithm for computing observables and to demonstrate the gain in accuracy over standard supervised deep learning algorithms, such as the one proposed in \cite{LMR1} . We will also use this algorithm in the context of speeding up forward uncertainty quantification. Numerical experiments for two prototypical problems in scientific computing will be presented in order to illustrate the gain in efficiency with the proposed algorithm.

The rest of the paper is organized as follows: in section \ref{sec:2}, we present the deep learning algorithm for approximating observables. The multi-level machine learning algorithm is presented in section \ref{sec:3} and its application to uncertainty quantification (UQ) is presented in section \ref{sec:4}. In section \ref{sec:5}, we present extensions of the multi-level algorithm and discuss its implementation. Numerical experiments are presented in section \ref{sec:6}. 

\section{The deep learning algorithm}
\label{sec:2}
\subsection{Problem formulation}
\label{sec:21}
Our objective is to approximate observables with machine learning algorithms. For definiteness, we assume that the observable of interest is defined in terms of the solutions of the following very generic system of time-dependent parametric PDEs:
\begin{equation}
\label{eq:ppde}
\begin{aligned}
\partial_t \U(t,x,y) &= L\left(y, \U, \nabla_x \U, \nabla^2_x \U, \ldots \right), \quad \forall~(t,x,y) \in [0,T] \times D(y) \times Y, \\
\U(0,x,y) &= \overline{\U}(x,y), \quad \forall~ (x,y) \in D(y) \times Y, \\
L_{b} \U(t,x,y) &= \U_b (t,x,y), \quad \forall ~(t,x,y) \in [0,T] \times \partial D(y) \times Y.
\end{aligned}
\end{equation}
Here, $Y$ is the underlying parameter space and without loss of generality, we assume it to be $Y = [0,1]^{d}$, for some $d \in {\mathbb N}$. 

The spatial domain is labeled as $y \rightarrow D(y) \subset \R^{d_s}$ and  $\U: [0,T] \times D \times Y  \rightarrow \R^m$ is the vector of unknowns. The differential operator $L$ is in a very generic form and can depend on the gradient and Hessian of $\U$, and possibly higher-order spatial derivatives. For instance, the heat equation as well as the Euler or Navier-Stokes equations of fluid dynamics are specific examples. Moreover, $L_b$ is a generic operator for imposing boundary conditions. 

The parametric nature of the PDE \eqref{eq:ppde}, represented by the parameter space $Y$, can stem from uncertainty quantification or Bayesian inversion problems where the parameter space models uncertainties in the PDE. The parametric nature can also arise from optimal design, control and PDE constrained optimization problems with $Y$ being the design (control) space. We equip this parameter space with a measure $\mu \in {\rm Prob} (Y)$.  

For the parameterized PDE \eqref{eq:ppde}, we consider the following generic form of observables,
\begin{equation}
\label{eq:obsp}
L_g(y,\U) := \int\limits_0^T\int\limits_{D_y} \psi(x,t) g(\U(t,x,y)) dx dt, \quad {\rm for} ~\mu~{\rm a.e}~y \in Y.
\end{equation} 
Here, $\psi \in L^1_{{\rm loc}} (D_y \times (0,T))$ is a  \emph{test function} and $g \in C^s(\R^m)$, for $s \geq 1$. 

For fixed functions $\psi,g$, we define the \emph{parameters to observable} map:
\begin{equation}
\label{eq:ptoob}
\map:y \in Y \rightarrow \map(y) = L_g(y,\U),
\end{equation}
with $L_g$ being defined by \eqref{eq:obsp}. 

We also assume that there exist suitable numerical schemes for approximating the PDE \eqref{eq:ppde} for every parameter vector $y \in Y$, such that a high-resolution approximate solution $\U^{\Delta}(y) \approx \U(y)$ is available, with $\Delta$ denoting the grid resolution (mesh size, time step etc.). Hence, there exists an approximation to the \emph{inputs to observable} map $\map^\Delta$ of the form,
 \begin{equation}
\label{eq:ptoob1}
\map^\Delta:y \in Y \rightarrow \map^\Delta(y) = L_g(y,\U^\Delta),
\end{equation}
with the integrals in \eqref{eq:obsp} being approximated to high accuracy by quadratures. Therefore, the original input parameters to observable map $\map$ is approximated by $\map^{\Delta}$ to very high accuracy i.e, for every value of a tolerance $\epsilon > 0$, there exists a $\Delta \ll 1$, such that
\begin{equation}
\label{eq:err1}
\|\map(y) - \map^{\Delta}(y)\|_{L_{\mu}^p(Y)} < \epsilon,
\end{equation}
for some $1 \leq p \leq \infty$ and weighted norm, $$
\|f\|_{L_{\mu}^p(Y)} := \left(\int_Y |f(y)|^p d\mu(y) \right)^{\frac{1}{p}}.
$$
The estimate \eqref{eq:err1} can be ensured by choosing $\Delta$ small enough in an underlying error estimate,
\begin{equation}
\label{eq:err2}
\|\map(y) - \map^{\Delta}(y)\|_{L_{\mu}^p(Y)} \sim \Delta^s, 
\end{equation}
for some $s > 0$.
\subsection{Deep learning the parameters to observable map}
\label{sec:22}
The process of learning the (approximate) parameters to observable map $\map^{\Delta}$ \eqref{eq:ptoob1} involves the following steps:
 \subsubsection{Training set.} 
 \label{sec:training_set}
 As is customary in supervised learning (\cite{DLbook} and references therein), we need to generate or obtain data to train the network. To this end, we fix $N \in \N$ and select a set of parameters $\EuScript{S} = \{y_i\}_{1 \leq i \leq N}$, with each $y_i \in Y$. The points in $\EuScript{S}$ can be chosen randomly from the parameter space $Y$, independently and identically distributed with the measure $\mu$. We will identify the training set $\train$ with the vector $\vtrain \in Y^N$, defined by 
 $$
 \vtrain = \left[y_1,y_2,\ldots y_N \right].
 $$
 Hence, each $\vtrain \subset Y$, is distributed according to the measure $\mu_N \in \Prob(Y^N)$ with $\mu_N(y_1,y_2,\ldots,y_N)=\mu(y_1)\otimes\mu(y_2)\ldots\otimes\mu(y_N)$ 
 
 Once the training set $\train$ is chosen, we perform a set of simulations of the underlying PDE \eqref{eq:ppde}, at a resolution $\Delta$, to obtain $\map^{\Delta} (y)$, for all $y \in \train$.
 \begin{figure}[htbp]
\centering
\includegraphics[width=8cm]{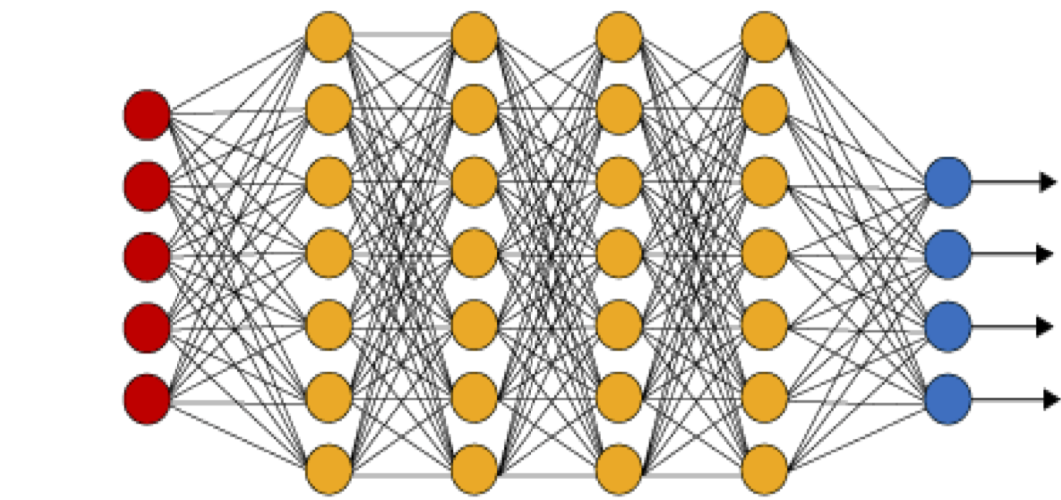}
\caption{An illustration of a (fully connected) deep neural network. The red neurons represent the inputs to the network and the blue neurons denote the output layer. They are
connected by hidden layers with yellow neurons. Each hidden unit (neuron) is connected by affine linear maps between units in different layers and then with nonlinear (scalar) activation functions within units.}
\label{fig:1}
\end{figure}

\subsubsection{Neural network.} 
\label{sec:NN}
Given an input vector $y \in Y$, a feedforward neural network (also termed as a multi-layer perceptron), shown in figure \ref{fig:1}, transforms it to an output, through a layer of units (neurons) which compose of either affine-linear maps between units (in successive layers) or scalar non-linear activation functions within units \cite{DLbook}, resulting in the representation,
\begin{equation}
\label{eq:ann1}
\map_{\theta}(y) = C_K \circ\sigma \circ C_{K-1}\ldots \ldots \ldots \circ\sigma \circ C_2 \circ \sigma \circ C_1(y).
\end{equation} 
Here, $\circ$ refers to the composition of functions and $\sigma$ is a scalar (non-linear) activation function. A large variety of activation functions have been considered in the machine learning literature \cite{DLbook}. A very popular choice, which we will consider for the rest of this article, is the \emph{ReLU} function,
\begin{equation}
\label{eq:relu}
\sigma(z) = \max(z,0).
\end{equation}
When, $z \in \R^p$ for some $p > 1$, then the output of the ReLU function in \eqref{eq:relu} is evaluated componentwise. 

For any $1 \leq k \leq K$, we define
\begin{equation}
\label{eq:C}
C_k z_k = W_k z_k + b_k, \quad {\rm for} ~ W_k \in \R^{d_{k+1} \times d_k}, z_k \in \R^{d_k}, b_k \in \R^{d_{k+1}}.
\end{equation}
For consistency of notation, we set $d_1 = d$ and $d_K = 1$. 

Thus in the terminology of machine learning (see also figure \ref{fig:1}), our neural network \eqref{eq:ann1} consists of an input layer, an output layer and $(K-1)$ hidden layers for some $1 < K \in \N$. The $k$-th hidden layer (with $d_k$ neurons) is given an input vector $z_k \in \R^{d_k}$ and transforms it first by an affine linear map $C_k$ \eqref{eq:C} and then by a ReLU (or another) nonlinear (component wise) activation $\sigma$ \eqref{eq:relu}. A straightforward addition shows that our network contains $\left(d + 1 + \sum\limits_{k=2}^{K-1} d_k\right)$ neurons. 
We also denote, 
\begin{equation}
\label{eq:theta}
\theta = \{W_k, b_k\},~ \theta_W = \{ W_k \},\quad \forall~ 1 \leq k \leq K,
\end{equation} 
to be the concatenated set of (tunable) weights for our network. It is straightforward to check that $\theta \in \Theta \subset \R^M$ with
\begin{equation}
\label{eq:ns}
M = \sum\limits_{k=1}^{K-1} (d_k +1) d_{k+1}.
\end{equation}

\subsubsection{Loss functions and optimization.} 
For any $y \in \train$, one can readily compute the output of the neural network $\map_{\theta} (y)$ for any weight vector $\theta \in \Theta$. We define the so-called training \emph{loss function} as 
\begin{equation}
\label{eq:lf1}
J (\theta) : = \sum\limits_{y \in \train} |\map^{\Delta}(y) - \map_{\theta} (y) |^p,
\end{equation}
for some $1 \leq p < \infty$.  

The goal of the training process in machine learning is to find the weight vector $\theta \in \Theta$, for which the loss function \eqref{eq:lf1} is minimized. 

It is common in machine learning \cite{DLbook} to regularize the minimization problem for the loss function i.e we seek to find,
\begin{equation}
\label{eq:lf2}
\theta^{\ast} = {\rm arg}\min\limits_{\theta \in \Theta} \left(J(\theta) + \lambda \reg(\theta) \right).
\end{equation}  
Here, $\reg:\Theta \to \R$ is a \emph{regularization} (penalization) term. A popular choice is to set  $\reg(\theta) = \|\theta_W\|^q_q$ for either $q=1$ (to induce sparsity) or $q=2$. The parameter $0 \leq \lambda \ll 1$ balances the regularization term with the actual loss $J$ \eqref{eq:lf1}. 

The above minimization problem amounts to finding a minimum of a possibly non-convex function over a subset of $\R^M$ for possibly very large $M$. We follow standard practice in machine learning by either (approximately) solving \eqref{eq:lf2} with a full-batch gradient descent algorithm or variants of mini-batch stochastic gradient descent (SGD) algorithms such as ADAM \cite{ADAM}. 

For notational simplicity, we denote the (approximate, local) minimum weight vector in \eqref{eq:lf2} as $\theta^{\ast}$ and the underlying deep neural network $\map^{\ast}= \map_{\theta^{\ast}}$ will be our neural network surrogate for the parameters to observable map $\map$ \eqref{eq:ptoob1}. The algorithm for computing this neural network is summarized below,
\begin{algorithm} 
\label{alg:DL} {\bf Deep learning of parameters to observable map}. 
\begin{itemize}
\item [{\bf Inputs}:] Parameterized PDE \eqref{eq:ppde}, Observable \eqref{eq:obsp}, high-resolution numerical method for solving \eqref{eq:ppde} and calculating \eqref{eq:obsp}.
\item [{\bf Goal}:] Find neural network $\map_{\theta^{\ast}}$ for approximating the parameters to observable map $\map$ \eqref{eq:ptoob1}. 
\item [{\bf Step $1$}:] Choose the training set $\train$ and evaluate $\map^{\Delta}(y)$ for all $y \in \train$ by a numerical method. 
\item [{\bf Step $2$}:] For an initial value of the weight vector $\overline{\theta} \in \Theta$, evaluate the neural network $\map_{\overline{\theta}}$ \eqref{eq:ann1}, the loss function \eqref{eq:lf2} and its gradients to initialize the
(stochastic) gradient descent algorithm.
\item [{\bf Step $3$}:] Run a stochastic gradient descent algorithm till an approximate local minimum $\theta^{\ast}$ of \eqref{eq:lf2} is reached. The map $\map^{\ast} = \map_{\theta^{\ast}}$ is the desired neural network approximating the
parameters to observable map $\map$.

\end{itemize}
\end{algorithm}
Note that the \emph{trained neural network} $\map^{\ast} = \map^{\ast}(\vtrain)$ depends explicitly on the training set, identified by the vector $\vtrain \in Y^N$. Hence, it should be denoted as $\map^{\ast}(y;\vtrain)$ for its application to each $y \in Y$. However for notational convenience, we will suppress this explicit dependence and denote the trained network as $\map^{\ast}$.
\subsection{An estimate on the generalization error of deep learning algorithm \ref{alg:DL}}
\label{sec:23}
For the rest of this section, we set $p=1$ in \eqref{eq:lf1} and aim to minimize the absolute value of the training loss. Moreover, we assume that their exists an underlying complete probability space $(\Omega, \Sigma, \IP)$, with respect to which random draws can be made.

Our aim in this section is to derive bounds on the so-called \emph{generalization error} \cite{MLbook}, which is customarily defined by,
\begin{equation}
    \label{eq:egen}
    \er_G(\vtrain) = \er_{G} (\theta^{\ast};\vtrain) := \int\limits_{Y} |\map^{\Delta}(y) - \map^{\ast}(y;\vtrain)| d\mu(y)
\end{equation}
Note that this generalization error depends explicitly on the training set $\train$ (identified by vector $\vtrain$). 

For each fixed training set $\train$, the training process in algorithm \ref{alg:DL} amounts to minimizing the so-called \emph{training error}:

\begin{equation}
\label{eq:etrain}
\er_T(\vtrain) = \er_{T} (\theta^{\ast};\vtrain) : = \frac{1}{N}\sum\limits_{i=1}^N |\map^{\Delta}(y_i) - \map^{\ast}(y_i;\vtrain) |,
\end{equation}
with $y_i \in \train$. The training error can be estimated from the loss function \eqref{eq:lf1}, a posteriori. Note that the training error $\er_T$ depends on the underlying randomly drawn training vector $\vtrain$. 

It would be tempting to realize that for each randomly drawn $\vtrain$, the training error \eqref{eq:etrain} in the collocation of the integrand $|\map^{\Delta} - \map^{\ast}|$, of the generalization error \eqref{eq:egen} on randomly chosen points. A naive application of the central limit theorem would then lead to a bound on the so-called \emph{generalization gap} of the form, 
\begin{equation}
    \label{eq:bn1}
    \int_{Y^N} |\er_G(\vtrain) - \er_T(\vtrain)|^2 d\mu^N(\vtrain) \sim \frac{\var\left(|\map^{\Delta} - \map^{\ast}(:;\vtrain)|\right)}{N}.
\end{equation}
Thus, in heuristic terms, the generalization gap is estimated in terms of the variance of the integrand in \eqref{eq:egen} and the number of training samples. 

However, such a bound is not rigorous. This is on account of the fact that the central limit theorem is only applicable if the underlying realizations are independent. Although this is true for the underlying map $\map^{\Delta}$, this is no longer true for the \emph{trained neural network} $\map^{\ast}$. In fact, the realizations of $\map^{\ast}$ on the training points can be  \emph{highly correlated} during the training process. These correlations create a formidable obstacle for obtaining sharp bounds on the generalization error \cite{MLbook}. In fact, tools from statistical learning theory (\cite{CS1}) such as VC dimension or Rademacher complexity \cite{MLbook} have been developed to deal with this issue, see the recent paper \cite{EMW1} for an application of this theory to obtain sharp generalization bounds for one hidden layer neural networks. 

Our objective in this paper is not to derive or work with sharp bounds on the generalization error but rather to illustrate the role of variance in the generalization error and how reducing variance can increase accuracy. To this end, we will adopt a heuristic approach and assume certain properties of trained neural networks that allow us to sharply illustrate the role of the underlying variance in controlling the generalization error.

To this end, we start by realizing that the generalization error \eqref{eq:egen} relies explicitly on the training set $\train$ and by defining an average (over all training sets) \emph{cumulative generalization error}: 
\begin{equation}
    \label{eq:ecgen}
    \bar{\er}_G = \int\limits_{Y^N} \er_{G} (\vtrain) d\mu^N(\vtrain) = \int\limits_{Y^N} \int\limits_{Y} |\map^{\Delta}(y) - \map^{\ast}(y;\vtrain)| d\mu(y)d\mu^N(\vtrain).
    \end{equation}
Similarly an average \emph{cumulative training error} is defined by,
\begin{equation}
\label{eq:ectrain}
\bar{\er}_T = \int_{Y^n} \er_{T}(\vtrain) d\mu^N(\vtrain) = \frac{1}{N}\int_{Y^N}\sum\limits_{i=1}^N |\map^{\Delta}(y_i) - \map^{\ast}(y_i;\vtrain) | d\mu^N(\vtrain).
\end{equation}
Note that as the points $y_i \in \train$ in the integrand of \eqref{eq:ectrain}, the cumulative training error is a deterministic quantity. 

Our objective would be to estimate the \emph{cumulative generalization gap} i.e the difference between \eqref{eq:ecgen} and \eqref{eq:ectrain}. To do so, we require another widely used set in machine learning i.e the so-called \emph{validation set},
\begin{equation}
    \label{eq:val}
    \val = \{z_j \in Y, ~ 1\leq j \leq N, \quad z_j ~i.i.d~wrt~\mu\}.
\end{equation}
The validation set is chosen before the starting of the training process and is independent of the training sets. We define the \emph{cumulative validation error} as, 
\begin{equation}
\label{eq:ecval}
\bar{\er}_V= \frac{1}{N}\int_{Y^N}\sum\limits_{j=1}^N |\map^{\Delta}(z_j) - \map^{\ast}(z_j;\vtrain) | d\mu^N(\vtrain).
\end{equation}
We observe that the as the set $\val$ is drawn randomly from $Y$ with underlying distribution $\mu$, the cumulative validation error is a random quantity, $\bar{\er}_V = \bar{\er}_V(\omega)$ with $\omega \in \Omega$. We suppress this $\omega$-dependence for notational convenience. Finally, we introduce the \emph{validation gap}:
\begin{equation}
    \label{eq:vgap}
    \er_{TV}:= \IE\left(|\bar{\er}_{T} - \bar{\er}_V|\right):= \int\limits_{\Omega}|\bar{\er}_{T} - \bar{\er}_V(\omega)|d\IP(\omega) 
\end{equation}
Equipped with the above notation and considerations, we obtain the following bounds on the cumulative generalization error,
\begin{lemma}
\label{theo:1}
The generalization gap i.e $\bar{\er}_G - \bar{\er}_T$, for the deep learning algorithm \ref{alg:DL} for approximating the observable $\map^{\Delta}$ \eqref{eq:ptoob1} satisfies the following bound,
\begin{equation}
    \label{eq:genE}
|\bar{\er}_G - \bar{\er}_T| \leq \er_{TV} +  
    \sqrt{\frac{8\left(\var\left(\map^{\Delta}\right)+\var\left(\map^{\ast}\right)\right)}{N}}.
\end{equation}
Here, $\var$ denotes the following variances,
\begin{equation}
\label{eq:var}
\begin{aligned}
\var\left(\map^{\Delta} \right) &= \int\limits_Y \left(\map^{\Delta}(y)\right)^2 d\mu(y) - \left( \int\limits_Y \map^{\Delta}(y) d\mu(y) \right)^2, \\
\var\left(\map^{\ast} \right) &= \int\limits_Y \int\limits_{Y^N} \left(\map^{\ast}(y;\vtrain)\right)^2 d\mu^N(\vtrain) d\mu(y) - \left( \int\limits_Y \int\limits_{Y^N} \map^{\ast}(y;\vtrain) d\mu^N(\vtrain) d\mu(y) \right)^2,
\end{aligned}
\end{equation}
\end{lemma}
\begin{proof}
We start with the elementary application of triangle and Cauchy-Schwartz inequalities to obtain,
\begin{align*}
    |\bar{\er}_G - \bar{\er}_T| = \IE(|\bar{\er}_G - \bar{\er}_T|) &\leq \er_{TV} + \IE(|\bar{\er}_G - \bar{\er}_V|)  \\ 
    &\leq \er_{TV} + \sqrt{\IE(|\bar{\er}_G - \bar{\er}_V|^2)}
\end{align*}

Next, we will use a Monte Carlo approximation to compute the integral inside the square root of the above expression. 

Using the definitions of the generalization error \eqref{eq:ecgen} and validation error \eqref{eq:ecval} and for any $z \in Y$, consider the integrand in the above integral i.e, $I(z) = \map^{\Delta}(z) - \int_{Y^N}\map^{\ast}(z;\vtrain) d\mu^N(\vtrain)$. As the validation points $z_j \in \val$ are randomly chosen from the underlying measure $\mu$, we realize that the validation error $\er_V$ \eqref{eq:ecval} is the Monte Carlo quadrature approximation of the generalization error \eqref{eq:ecgen} and we can estimate the difference between them, in terms of the central limit theorem \cite{CAF1} by,
\begin{equation*}
\IE(|\bar{\er}_G - \bar{\er}_V|^2) \leq \frac{\var(I)}{N}
\end{equation*}

We estimate the variance in the above equation in the following manner,
\begin{align*}
    \var(I) &= \int\limits_{Y} \left( \int_{Y}\int_{Y^N}|\map^{\Delta}(z) - \map^{\ast}(z;\vtrain)| - |\map^{\Delta}(\bar{z}) - \map^{\ast}(\bar{z};\vtrain)| d\mu^N(\vtrain) d\mu(\bar{z}) \right)^2 d\mu(z) \\
    &\leq \int\limits_Y \left(\int\limits_Y \int\limits_{Y^N}|\map^{\Delta}(z) - \map^{\ast}(z;\vtrain) - (\map^{\Delta}(\bar{z}) - \map^{\ast}(\bar{z};\vtrain))|d\mu^N(\vtrain) d\mu(\bar{z}) \right)^2 d\mu(z) ~({\rm by~triangle~inequality})\\
    &\leq \int\limits_Y \left(\int\limits_Y  |\map^{\Delta}(z) - \map^{\Delta}(\bar{z})| d\mu(\bar{z}) + \int\limits_Y \int\limits_{Y^N} |\map^{\ast}(z;\vtrain) - \map^{\ast}(\bar{z};\vtrain)| d\mu^N(\vtrain) d\mu(\bar{z}) \right)^2 d\mu(z)~({\rm by~triangle~inequality})\\ 
    &\leq 2\left( \underbrace{\int\limits_Y \int\limits_Y |\map^{\Delta}(z) - \map^{\Delta}(\bar{z})|^2 d\mu(z) d\mu(\bar{z})}_{I_1} + \underbrace{\int\limits_Y \int\limits_Y |\map^{\ast}(z;\vtrain) - \map^{\ast}(\bar{z};\vtrain)|^2 d\mu^N(\vtrain) d\mu(\bar{z}) d\mu(z)}_{I_2}\right)~({\rm by~Cauchy-Schwartz})\
\end{align*}
Using the definitions of the variances \eqref{eq:var}, it is straightforward to obtain the following estimates,
$$
I_1 \leq 4 \var\left(\map^{\Delta}\right), ~ I_2 \leq 4 \var\left(\map^{\ast}\right),
$$
resulting in the estimate \eqref{eq:genE}.
 \\
\end{proof}
The bound \eqref{eq:genE} yields
\begin{equation}
    \label{eq:genE2}
    \bar{\er}_G \sim \bar{\er}_T + \er_{TV} + 2\sqrt{2} \frac{\left({\rm std}(\map^{\Delta}) + {\rm std}(\map^{\ast})\right)}{\sqrt{N}}, \\
\end{equation}
with std denoting the standard deviation i.e, the square roof of the variances in  \eqref{eq:var}.

Several remarks are in order about the estimate \eqref{eq:genE} on the generalization gap and the resulting estimate \eqref{eq:genE2} on the generalization error. 
\begin{remark}
The estimate \eqref{eq:genE2} bounds the \emph{cumulative generalization error} (generalization error averaged over different choices of the training sets) in terms of the \emph{cumulative training error} (training error averaged over different choices of the training sets), the \emph{validation gap} (difference between training and validation errors, averaged over training sets), and the variances of the underlying map as well as the trained neural network (averaged over the choice of training samples). This is a non-standard form of writing the generalization error. However, it corresponds to the usual practice in machine learning, where one computes both the training and validation errors for each choice of training samples and also considers multiple training samples.
\end{remark}
\begin{remark}
In practice, the accuracy of the training process is ascertained by monitoring both the training and the validation errors. The training process is usually terminated when the training error is lower than some tolerance, while at the same time the validation error is also low enough, i.e the validation gap $\er_{TV}$ is small. So, an estimate of the form \eqref{eq:genE2} corresponds to what is most often computed in practice. One of the limitations of the estimate \eqref{eq:genE2} is that fact that the training and validation sets are of similar size whereas in practice, one sets aside a much smaller number of samples for validation. 
\end{remark}

\begin{remark}
\label{re:sharpness}
The above estimate \eqref{eq:genE2} on the generalization gap is clearly an upper bound and is not necessarily sharp. In fact, the proof of the lemma makes it clear that the bound could be a significant overestimate. As is well known \cite{ARORA,EMW1,NEYS1}, estimating the generalization error of neural networks sharply is a notoriously hard problem. Nevertheless, the upper bound \eqref{eq:genE2} will suffice for our purposes in this paper. 
\end{remark}

\section{A multi-level Deep learning algorithm.}
\label{sec:3}
In order to motivate the design of a multi-level deep learning algorithm, we introduce the following concept,
\begin{definition}
\label{def:wtnn}
{\bf Well-trained Neural Network.} A neural network $\map^{\ast}$, generated by the deep learning algorithm \ref{alg:DL} to approximate the parameters to observable map $\map^{\Delta}$ \eqref{eq:ptoob1}, is said to be well-trained if the following hold,
\begin{equation}
    \label{eq:desp1}
    {\rm std}(\map^{\ast}) \sim {\rm std}(\map^{\Delta}), \quad \bar{\er}_T \approx \er_{TV} \ll \frac{{\rm std}(\map^{\Delta})}{\sqrt{N}}.
\end{equation}
\end{definition}
In other words, the training error \eqref{eq:ectrain} and the validation gap \eqref{eq:vgap} for a \emph{well-trained network} are significantly lower than the variance of the underlying map. Moreover, the variance of the network is comparable to that of the underlying map. 

We assume that such an well-trained network can be found during the training process. Although this assumption appears rather stringent, it must be mentioned that the conditions on the training error and the validation gap can be monitored during the training process. In fact, it is standard practice in machine learning to check the validation error in order to monitor overfitting. 

This assumption automatically implies from \eqref{eq:genE2} that the generalization error scales as,
\begin{equation}
\label{eq:genE3}
\bar{\er}_G \sim \frac{{\rm std}(\map^{\Delta})}{\sqrt{N}}.
\end{equation}

Even under the assumption of the existence of an well-trained network, as long as ${\rm std}(\map^{\Delta}) \sim {\mathcal O}(1)$, we see from \eqref{eq:genE3} that the generalization error $\bar{\er}_G \sim \frac{1}{\sqrt{N}}$. Hence, for obtaining a generalization error of say $1$ percent, we need $10^4$ training samples. Generating such a large number of training samples might be prohibitively expensive, particularly for if the underlying parametric PDE \eqref{eq:ppde} is in two or three space dimensions. 

\begin{figure}[htbp]
\centering
\includegraphics[width=15cm]{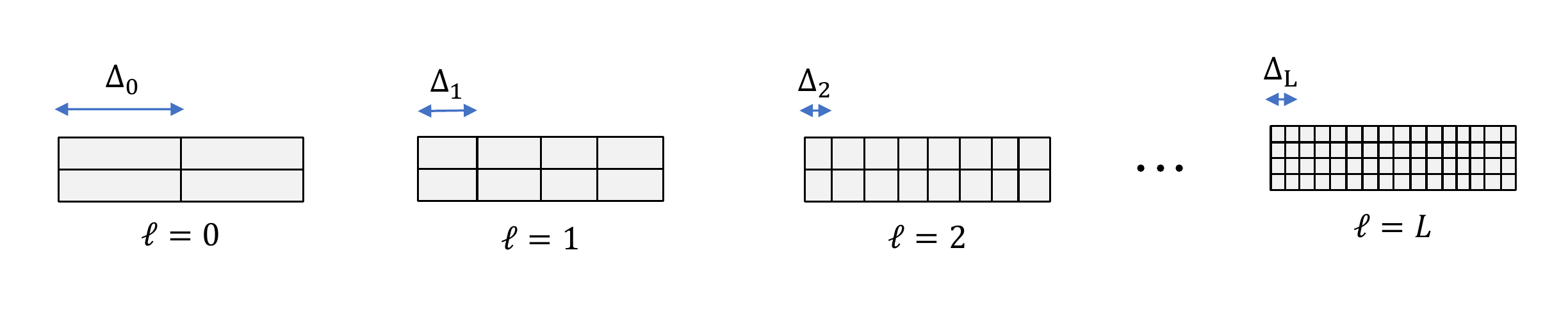}
\caption{A schematic for a sequence of nested grids used in defining the multi-level algorithms.}
\label{fig:2}
\end{figure}
Our goal in this section is to propose a \emph{multi-level} version of the deep learning algorithm \ref{alg:DL}. The basis of this algorithm is the observation that the underlying parameters to observable map $\map$ \eqref{eq:ptoob} can be approximated on a sequence of mesh resolutions $\Delta_\ell$, for $0 \leq \ell \leq L$ for some $L > 0$. We require that $\Delta_\ell < \Delta_{\ell -1}$ for each $1 \leq \ell \leq L$, see \figurename{ \ref{fig:2}} for a diagrammatic representation of this sequence of grids.

For each $\ell$, we assume that the underlying parameters to observable map $\map$ can be approximated by the map $\map^{\Delta_\ell}$, computed on resolution $\Delta_\ell$ with an estimate,
\begin{equation}
    \label{eq:app1}
    \|\map - \map^{\Delta_\ell}\|_{L^1_\mu(Y)} \sim \Delta_\ell^s, \quad s > 0.
\end{equation}

Given such a sequence of resolutions $\Delta_\ell$ and approximate parameters to observable maps $\map^{\Delta_\ell}$, we have the following (telescopic decomposition)
\begin{equation}
    \label{eq:tsum1}
    \map^{\Delta_L}(y) = \map^{\Delta_0}(y) + \sum\limits_{\ell=1}^L\left( \map^{\Delta_\ell}(y) - \map^{\Delta_{\ell-1}}(y)\right), \quad \forall y \in Y. 
    \end{equation}
Introducing the so-called \emph{details},
\begin{equation}
    \label{eq:dtl1}
    \dtl_\ell(y) =  \map^{\Delta_\ell}(y) - \map^{\Delta_{\ell-1}}(y), \quad \forall y 
    \in Y,~ 1 \leq \ell \leq L,
\end{equation}
we rewrite the telescopic decomposition \eqref{eq:tsum1} as
\begin{equation}
    \label{eq:tsum}
    \map^{\Delta_L}(y) = \map^{\Delta_0}(y) + \sum\limits_{\ell=1}^L \dtl_\ell(y), \quad \forall y \in Y. 
    \end{equation}

The multi-level deep learning algorithm will be based on independently learning the following maps \begin{equation}
    \label{eq:ml1}
    \map^{\Delta_0} (y) \approx \map_0^{\ast} (y), \quad \dtl_\ell (y) \approx \dtl^{\ast}_{\ell}(y), \quad \forall y\in Y, ~ 1 \leq \ell \leq L.
\end{equation}    
Here, $\map_0^{\ast},\dtl^{\ast}_\ell$, for $1 \leq \ell \leq L$ are artificial neural networks of the form \eqref{eq:ann1}. The resulting algorithm is summarized below,

\begin{algorithm} 
\label{alg:ML} {\bf Multi-level deep learning of parameters to observable map}. 
\begin{itemize}
\item [{\bf Inputs}:] Parameterized PDE \eqref{eq:ppde}, Observable \eqref{eq:obsp}, high-resolution numerical method for solving \eqref{eq:ppde} and calculating \eqref{eq:obsp}, a sequence of grids with grid size $\Delta_\ell$ for $0 \leq \ell \leq L$.
\item [{\bf Goal}:] Find neural network $\map^{\ast}_{ml}$ for approximating the parameters to observable map $\map$ \eqref{eq:ptoob}. 
\item [{\bf Step $1$}:] For the coarsest mesh resolution $\Delta_0$, select a training set $\train_0 = \{y^0_i\},~1 \leq i \leq N_0 = \#(\train_0)$, with $y^0_i \in Y$ being independently and identically distributed with respect to the measure $\mu$ i.e, the training set $\train_0$ can be identified with the $\vtrain_0 \in Y^{N_0}$, with 
$$
\vtrain_0 = \left[y^0_1,y^0_2,\ldots,y^0_N\right],
$$
drawn from an underlying distribution $\mu^{N_0} \in \Prob(Y^{N_0})$. For each $y^0_i$, compute $\map^{\Delta}_0(y^0_i)$ by solving the PDE \eqref{eq:ppde} on a mesh resolution $\Delta_0$ and computing the observable \eqref{eq:obsp}. With the training data $\{y^0_i,\map^{\Delta_0}(y^0_i)\}_{1 \leq i \leq N_0}$, find the neural network $\map^{\ast}_0 \approx \map^{\Delta}_0$ by applying the deep learning algorithm \ref{alg:DL}. Note that the trained neural network $\map^{\ast}_0 = \map^{\ast}_0(:;\vtrain_0)$, depends on the underlying training set and we suppress this dependence for notational simplicity. 
\item [{\bf Step $2$}:] For each $1 \leq \ell \leq L$, select a training set $\train_\ell = \{y^\ell_i\},~1 \leq i \leq N_\ell = \#(\train_\ell)$, with $y^\ell_i \in Y$ being independently and identically distributed with respect to the measure $\mu$, i.e, the training set $\train_\ell$ can be identified with the vector $\vtrain_\ell \in Y^{N_\ell}$ as
$$
\vtrain_\ell = \left[y^\ell_1,y^\ell_2,\ldots,y^\ell_N\right],
$$

drawn from an underlying distribution $\mu^{N_\ell} \in \Prob(Y^{N_\ell})$. For each $y^\ell_i$, compute $\dtl_\ell(y^\ell_i)= \map^{\Delta_\ell}(y^\ell_i) - \map^{\Delta_{\ell-1}}(y^\ell_i)$ by solving the PDE \eqref{eq:ppde} on two mesh resolutions of mesh size $\Delta_\ell$ and $\Delta_{\ell-1}$ and computing the observable \eqref{eq:obsp}. With the training data $\{y^\ell_i,\dtl_\ell(y^\ell_i)\}_{1 \leq i \leq N_\ell}$, find the neural network $\dtl^{\ast}_\ell \approx \dtl_\ell$ by applying the deep learning algorithm \ref{alg:DL}. Note that the trained neural network $\dtl^{\ast}_\ell = \dtl^{\ast}_\ell(:;\vtrain_\ell)$, depends on the underlying training set and we suppress this dependence for notational simplicity.
\item [{\bf Step $3$}:] Form the artificial neural network $\map^{\ast}_{ml} \approx \map$ by,
\begin{equation}
    \label{eq:ml}
    \map_{ml}^{\ast}(y) = \map^{\ast}_0(y) + \sum\limits_{\ell=1}^L \dtl^\ast_\ell(y), \quad \forall y \in Y.
\end{equation}

\end{itemize}
\end{algorithm}
Note that the neural network $\map^{\ast}_{ml}$ explicitly depends on the underlying training sets. To formalize this dependence, we introduce the notation,
\begin{align*}
\bar{N} &= N_0 + \sum\limits_{\ell=1}^L N_l, \quad Y^{\bar{N}} = Y^{N_0}\otimes Y^{N_1}\otimes \ldots \otimes Y^{N_L}, \\
\mu^{\bar{N}} &\in \Prob(Y^{\bar{N}}), \quad \mu^{\bar{N}} = \mu^{N_0}\otimes\mu^{N_1}\otimes \ldots \otimes \mu^{N_L}, \quad \bar{\vtrain} = \left[\vtrain_0,\vtrain_1,\ldots,\vtrain_L \right]
\end{align*}
With the above notation, we have that $\map^{\ast}_{ml} = \map_{ml}^{\ast}(:;\bar{\vtrain})$. 

Analogous to \eqref{eq:ecgen}, we can define a \emph{cumulative generalization error} for the multi-level network $\map^{\ast}_{ml}$ as,
\begin{equation}
\label{eq:genEml1}
{\bar{\er}}_G^{ml} := \int\limits_{Y^{\bar{N}}} \int\limits_Y |\map^{\Delta_L}(y) - \map^{\ast}_{ml}(y;\bar{\vtrain})| d\mu^{\bar{N}}(\bar{\vtrain}) d\mu(y).     
\end{equation}
Implicitly, we have assumed that the approximation $\map^{\Delta_L}$ to the underlying parameters to observable map $\map$ at the finest level of resolution $\Delta_L$ is the \emph{ground truth}. 

We are going to estimate the cumulative generalization error \eqref{eq:genEml1} in terms of the cumulative generalization errors at the coarsest level of resolution and that of the details, defined as
\begin{equation}
    \label{eq:genEml2}
    \bar{\er}_G^{0} = \int\limits_{Y^{N_0}} \int\limits_Y |\map^{\Delta_0}(y) - \map^{\ast}_{0}(y;\vtrain_0)| d\mu(y) d\mu^{N_0}(\vtrain_0), \quad  
    \bar{\er}_G^{\ell} = \int\limits_{Y^{N_\ell}} \int\limits_Y |\dtl_\ell(y) - \dtl^{\ast}_{\ell}(y;\vtrain_\ell)| d\mu(y) d\mu^{N_\ell}(\vtrain_\ell), \quad 1 \leq \ell \leq L.
\end{equation}
The cumulative training errors for each network can be defined analogous to \eqref{eq:ectrain} to obtain $\bar{\er}^0_T, \bar{\er}^\ell_T$. As in the previous section, we need to specify \emph{validation sets} of the form \eqref{eq:val}, for the validating the networks $\map_0^{\ast},\dtl^\ast_{\ell}$ for each $1 \leq \ell \leq L$. The \emph{validation gap} is defined, in analogy with \eqref{eq:ecval}, to define $\er_{TV}^\ell$, with $0 \leq \ell \leq L$.

We have the following bound for the cumulative generalization error \eqref{eq:genEml1} of the network, generated by the multi-level deep learning algorithm \ref{alg:ML}. 
\begin{lemma}
\label{lem:1} The generalization error \eqref{eq:genEml1} of the multi-level machine learning algorithm \ref{alg:ML} is estimated as,
\begin{equation}
    \label{eq:genEml}
    \begin{aligned}
    \bar{\er}_G^{ml} &\leq \bar{\er}^0_G + \sum\limits_{\ell=1}^L \bar{\er}^\ell_G, \\
     |\bar{\er}^0_G - \bar{\er}^0_T| &\leq \er^0_{TV} +  
    \sqrt{\frac{8\left(\var\left(\map^{\Delta_0}\right)+\var\left(\map^{\ast}_0\right)\right)}{N_0}}, \\
     |\bar{\er}^{\ell}_G - \bar{\er}^{\ell}_T| &\leq \er^{\ell}_{TV} +  
    \sqrt{\frac{8\left(\var\left(\dtl_{\ell}\right)+\var\left(\dtl^{\ast}_{\ell}\right)\right)}{N_\ell}}\quad {\rm for~all} \quad 1 \leq \ell \leq L.
    \end{aligned}
\end{equation}
\end{lemma}
The proof of the first inequality in \eqref{eq:genEml} in above lemma is based on the telescopic decomposition \eqref{eq:tsum} and successive applications of the triangle inequality. The other estimates in \eqref{eq:genEml} follow from a direct application of the inequality \eqref{eq:genE} in Lemma \ref{theo:1} to the neural networks $\map^{\ast}_0,\dtl^\ast_\ell$, for all $1 \leq \ell \leq L$. 

One can rewrite the above estimate on generalization error as,
\begin{equation}
    \label{eq:genEml3}
    \begin{aligned}
    \bar{\er}_G^{ml} &\sim \bar{\er}_T^0 + \er_{TV}^0 + \sum\limits_{\ell=1}^L \left(\bar{\er}_T^\ell + \er_{TV}^{\ell}\right) \\ 
    &+ \frac{2\sqrt{2}({\rm std}(\map^{\Delta_0}) + {\rm std}(\map^{\ast}_{0}))}{\sqrt{N_0}} + \sum\limits_{\ell=1}^L 
    \frac{2\sqrt{2}({\rm std}(\dtl_\ell) + {\rm std}(\dtl^{\ast}_\ell))}{\sqrt{N_l}}
    \end{aligned}
\end{equation}

We need to determine the training sample numbers $N_{\ell}$ for $0 \leq \ell \leq L$, in order to complete the description of the multi-level machine learning algorithm \ref{alg:ML}. These have to be determined such that the multi-level deep learning algorithm \ref{alg:ML} can lead to a greater accuracy (a lower generalization error $\bar{\er}_G$) at a similar computational cost as the underlying deep learning algorithm \ref{alg:DL}. 

A heuristic argument for selecting sample numbers runs as follows. First, we assume that the cost of training the neural networks and evaluating them is significantly lower than the cost of generating the training data with a PDE solver for \eqref{eq:ppde}. This assumption is indeed justified for most realistic problems (see table 13 in \cite{LMR1} for the training and evaluation costs vis a vis the cost of generating training data for a flow past airfoils). Next, we assume that the cost of solving a PDE such as \eqref{eq:ppde} for a single realization of the parameter vector $y \in Y$, on a mesh resolution of $\Delta_\ell$ scales as $\Delta_\ell^{-(d_s+1)}$. This assumption is justified for first-order time-dependent PDEs (due to the CFL condition) and the arguments below can be readily extended to a more general case. It assume it here for the sake of definiteness.

Thus, it is much cheaper to generate a training sample on a coarse resolution than on a fine resolution. Hence, the intuitive idea is to generate a much larger number of training samples on coarser mesh resolutions than on finer resolutions (see Figure \ref{fig:2}). From the generalization error estimate \eqref{eq:genEml3}, we see that contribution to the generalization error from the details, corresponding to fine mesh resolutions, can be very low, even for a low number of training samples as long as the standard deviation of the underlying details is low. Thus, we can combine a large number of training samples at low resolution with a few training samples at high resolution, in order to obtain low generalization errors at comparable cost to the deep learning algorithm \ref{alg:DL}.
These heuristic considerations are formalized in the lemma below.
\begin{lemma}
\label{lem:2}
For any given tolerance $\epsilon$, there exists a $\Delta > 0$ such that the error estimate \eqref{eq:err1} holds, we assume that training process for the deep learning algorithm \ref{alg:DL} results in a well-trained neural network $\map^{\ast}$ with properties \eqref{eq:desp1}. 
Moreover, we consider a generic sequence of mesh resolutions $s_n = \{\ell_k\}$ with $\ell_0=0$, $\ell_n=L$ that yields the following generalized formulation of the multi-level deep learning algorithm \ref{alg:ML}
\begin{equation}
\label{eq:def_mlml_gen}
\map_{\ell_n}^\ast(y) = \map_{\ell_0}^\ast(y) + \sum\limits_{k=1}^n \dtl^\ast_{k}(y),\quad \map_{\ell_0}^\ast(y)\approx \map^{\Delta_{0}}(y), \quad \dtl^{\ast}_{k}(y) \approx \dtl_{k} (y)  := \map^{\Delta_{\ell_k}}(y) - \map^{\Delta_{\ell_{k-1}}}(y),  \quad \forall y \in Y,
\end{equation}
such that
\begin{equation}
    \label{eq:sr1}
	\Delta_k = 2^{(\ell_n -\ell_k)} \Delta,\quad 0 \leq k \leq n,
\end{equation}
 and training sample numbers given by,
\begin{equation}
    \label{eq:tsn}
    N_k\sim \frac{L\var_k}{\epsilon^2}, \quad 0 \leq k \leq n, \quad {\rm with}~ \var_0 = \var(\map^{\Delta_0}),~\var_k = \var(\dtl_k), \quad 1 \leq k \leq n. 
\end{equation}
Furthermore, we assume that the training process in the multi-level deep learning algorithm \ref{alg:ML} results in \emph{well-trained} artificial neural networks i.e, neural networks $\map_0^{\ast},\dtl_k^{\ast}$ with properties,
\begin{equation}
    \label{eq:et2}
    \bar{\er}^0_T \approx \er_{TV}^0 \ll \frac{{\rm std}(\map^{\Delta_0})}{\sqrt{N_0}},~{\rm std}(\map^{\ast}_0) \sim {\rm std}(\map^{\Delta_0}),~\bar{\er}^k_T \approx \er_{TV}^k \ll \frac{{\rm std}(\dtl_k)}{\sqrt{N_k}},~ {\rm std}(\dtl^{\ast}_{k}) \sim {\rm std}(\dtl_k),  \forall 1 \leq k \leq n.
\end{equation}

Let $\Sigma_{ml}$ be speed up i.e, the ratio of computational cost of computing $\map^\Delta$ to accuracy of ${\mathcal O}(\epsilon)$. to the cost of computing $\map^\Delta$ to the same accuracy with the multi-level machine learning algorithm \ref{alg:ML}. Under the assumption that the costs of training and evaluation of all the neural networks in algorithms \ref{alg:DL} and \ref{alg:ML} is significantly smaller than the cost of generating the training data, we have the following estimate on $\Sigma_{ml}$,
\begin{equation}
    \label{eq:sp1}
    \frac{1}{\Sigma_{ml}} \sim \frac{L}{\var(L^{\Delta})} \left[  \var_0 2^{-L \bar{d}} + \sum\limits_{k =1}^n \var_k 2^{-(L-\ell_k) \bar{d}}    \right],
\end{equation}
with $\bar{d} = d_s +1$ is the number of space-time dimensions for the parametrized PDE \eqref{eq:ppde}. 
\end{lemma}
\begin{proof}
Based on error estimate \eqref{eq:err2}, we observe that $\epsilon \sim \Delta^s$. Applying  assumptions \eqref{eq:desp1} on the neural network $\map^{\ast}$, generated by the deep learning algorithm \ref{alg:DL}, in the estimate \eqref{eq:genE} for the generalization error of $\map^{\ast}$, we obtain that 
\begin{equation}
    \label{eq:N1}
    N \sim \var(\map^{\Delta}) \Delta^{-2s}
\end{equation}
training samples are required to approximate $\map^{\Delta}$ to tolerance $\epsilon$. 

Under our assumptions, the cost of training and evaluating neural networks is much smaller than the cost of generating the training data. Hence, the total cost of the deep learning algorithm \ref{alg:DL} is given by,
\begin{equation}
    \label{eq:cost1}
    \cost_{DL} \sim N \cost_{\Delta} \sim \var(\map^{\Delta}) \Delta^{-(2s+\bar{d})}.
\end{equation}
Here, we have used \eqref{eq:N1} and the fact that the computational cost of generating a single training sample by solving the PDE \eqref{eq:ppde} numerically at resolution $\Delta$ is given by $\cost_{\Delta} \sim \Delta^{-\bar{d}}$

Next, we consider the multi-level machine learning algorithm \ref{alg:ML}. Applying the assumptions \eqref{eq:sr1} on the resolutions, \eqref{eq:tsn} on the training sample numbers and \eqref{eq:et2} on trained neural networks in the formula for the generalization error \eqref{eq:genEml}, we observe that each term in the right hand side of the upper bound in \eqref{eq:genEml} scales as ${\mathcal O}\left(\frac{\epsilon}{n}\right)$ and the total generalization error is of ${\mathcal O}(\epsilon)$.

Under the assumptions that the cost of training and evaluating neural networks is much smaller than the cost of generating the training data, the total cost of the multi-level machine learning algorithm \ref{alg:ML} is given by,
\begin{equation}
    \label{eq:cost2}
    \begin{aligned}
    \cost_{ML} &\sim \sum\limits_{k=0}^n N_k \cost_{\Delta_{k}} \\
    &\sim \sum\limits_{k=0}^n L \var_k 2^{-(L-\ell_k) \bar{d}} \Delta^{-(2s+\bar{d})}, \quad {\rm by}~\eqref{eq:sr1},\eqref{eq:tsn}.
    \end{aligned}
\end{equation}
The estimate on speed-up \eqref{eq:sp1} follows from dividing \eqref{eq:cost2} by \eqref{eq:cost1}.

\end{proof}

\begin{remark}
Each sequence of resolutions in the multi-level method can be specified in terms of a single parameter, i.e \textit{model complexity} defined as:
\begin{equation}
\label{eq:cml}
    c_{ml}=\frac{n^2}{L}
\end{equation}
\end{remark}

\begin{remark}
From formula \eqref{eq:sp1}, we see that speedup of the multi-level algorithm \ref{alg:ML} over the underlying machine learning algorithm \ref{alg:DL} is expressed in terms of the \emph{variance of the details} $\dtl_{\ell}$. In practice, $\var_0 = \var(\map^{\Delta_0}) \sim \var(\map^{\Delta})$. On the other hand, a straightforward application of the error estimate \eqref{eq:err2} yields the following estimate,

\begin{equation}
    \label{eq:vard}
    \var_k = \var(\dtl_k) = \var(\map^{\Delta_{\ell_k}} - \map^{\Delta_{\ell_{k-1}}}) \sim \Delta^{2s}_{\ell_{k-1}} = \Delta^{2s} 2^{2(L - \ell_{k-1})s}, \quad {\rm from}~\eqref{eq:sr1}.
\end{equation}
Assuming for simplicity that $\ell_k = k\frac{L}{n}$ $\forall k: 1\leq k\leq n$ and substituting \eqref{eq:vard} in \eqref{eq:sp1} results in
\begin{equation}
    \label{eq:sp2}
     \frac{1}{\Sigma_{ml}} \sim
     L 2^{-L\bar{d}} + \frac{L\Delta^{2s}2^{2s\frac{L}{n}}}{\var(\map^{\Delta})}\sum\limits_{k=1}^n 2^{L(1-\frac{k}{n})(2s -\bar{d})}
\end{equation}
The geometric series in \eqref{eq:sp2} clearly converges to a finite value as long as $s \leq \frac{\bar{d}}{2}$. This holds true in most cases of practical interest as $s \leq 1$ and $\bar{d} = 3$ or $4$. Consequently, the speedup will be exponential in $L$ for small values of $\Delta$.
\end{remark}
\begin{remark}
We have used a very crude strategy of selecting sample numbers in \eqref{eq:tsn}. It relies on specifying the number of samples at the finest level $N_L$ and at the coarsest level $N_0$. The samples at intermediate levels are determined by,
\begin{equation}
    \label{eq:snumb}
    N_k = N_L2^{e(L-\ell_k)}, ~0 < k <n, \quad  e = \frac{\log_2(N_0/N_L)}{L}.
\end{equation}
More sophisticated strategies, such as those proposed in the context of multi-level Monte Carlo methods in \cite{GIL2} and references therein, might lead to greater speedup with the multi-level machine learning algorithms. 
\end{remark}

\section{Uncertainty Quantification}
\label{sec:4}
In this article, we focus our attention on forward uncertainty quantification (UQ) or uncertainty propagation with respect to the parameters to observables map $\map$ \eqref{eq:ptoob} (or rather its numerical surrogate $\map^{\Delta}$ \eqref{eq:ptoob1}). To this end, we follow \cite{LMR1} and consider the so-called \emph{push forward measure} with respect to this map, i.e, $\hat{\mu}^{\Delta} \in {\rm Prob}(\R)$ given by 
\begin{equation}
\label{eq:pf1} 
\hat{\mu}^{\Delta} := \map^{\Delta}\#\mu, \quad \Rightarrow \quad \int_{\R} f(z) d\hat{\mu}^{\Delta}(z) = \int_{Y} f(\map^{\Delta}(y)) d\mu(y), 
\end{equation}
for any $\mu$-measurable function $f:\R \to \R$. 

Note that the measure $\hat{\mu}^\Delta$ contains all the statistical information about the map $\map^{\Delta}$. In particular, any moment or statistical quantity of interest with respect to this map can be computed by choosing a suitable test function $f$ in \eqref{eq:pf1}. In particular, choosing $f(z) = z $ provides the mean $\bar{\map}$ of the observable $\map^{\Delta}$ and $f(z) = (z - \bar{\map})^2$ yields the variance.

The baseline Monte Carlo algorithm for approximating this measure (probability distribution) consists of choosing $J$ independent, identically distributed (with respect to $\mu$) samples $y_j \in Y$ and approximating the measure $\hat{\mu}^{\Delta}$ by the so-called \emph{empirical measure},
\begin{equation}
\label{eq:m1}
\hat{\mu}_{mc} = \frac{1}{J}\sum\limits_{j=1}^{J} \delta_{\map^{\Delta}(y_j)} \quad \Rightarrow \quad \int_{\R} f(z) d\hat{\mu}_{qmc}(z) = \frac{1}{J}\sum\limits_{j=1}^{J} f\left(\map^{\Delta}(y_j) \right).
\end{equation}

On the other hand, within the \emph{deep learning Monte Carlo} (DLMC) algorithm, proposed in \cite{LMR1}, we first generate a deep neural network $\map^{\ast}$ to approximate the underlying map $\map^{\Delta}$ by the deep learning algorithm \ref{alg:DL} and then approximate the measure $\hat{\mu}^{\Delta}$ by,
\begin{equation}
\label{eq:dlmc1}
\hat{\mu}^{\ast}_{mc} = \frac{1}{J_L}\sum\limits_{j=1}^{J_L} \delta_{\map^{\ast}(y_j)},
\end{equation}
with $J_L \gg J$ evaluations of the neural network. In \cite{LMR1}, the authors provide a complexity analysis of the DLMC algorithm and also estimate the speedup with respect to the baseline Monte Carlo algorithm, in terms of the errors of computing the measure  $\hat{\mu}^{\Delta}$ in the Wasserstein metric. We refer the interested reader to \cite{LMR1}, Theorem 3.11. The numerical experiments, presented in \cite{LMR1} also provide evidence of a significant speedup with the DLMC algorithm over the baseline Monte Carlo method. 

Here, we propose another variant of the DLMC algorithm, which is based on the multi-level Machine learning algorithm \ref{alg:ML}. The algorithm is as follows,
\begin{algorithm} 
\label{alg:ML2MC} {\bf A Multi-level Machine learning Monte Carlo (ML2MC) algorithm for forward UQ}. 
\begin{itemize}
\item [{\bf Inputs}:] Parameterized PDE \eqref{eq:ppde}, Observable \eqref{eq:obsp}, high-resolution numerical method for solving \eqref{eq:ppde} and calculating \eqref{eq:obsp}, a sequence of grids with grid size $\Delta_\ell$ for $0 \leq \ell \leq L$.
\item [{\bf Goal}:] Find a measure $\hat{\mu}_{ml2mc} \in {\rm Prob}(\R)$ to approximate the push-forward measure $\hat{\mu}^{\Delta}$ \eqref{eq:pf1}.
\item [{\bf Step $1$}:] Generate the neural network $\map^{\ast}_{ml} \approx \map^{\Delta}$ by applying the multi-level machine learning algorithm \ref{alg:ML}. 
\item [{\bf Step $2$}:] Define the approximate push-forward measure by,
\begin{equation}
\label{eq:ml2mc1}
\hat{\mu}^{\ast}_{ml2mc} = \frac{1}{J_L}\sum\limits_{j=1}^{J_L} \delta_{\map_{ml}^{\ast}(y_j)},
\end{equation}
\end{itemize}
\end{algorithm}
One can readily perform a complexity analysis, completely analogous to section 3.3 of \cite{LMR1} to quantify possible speedups with the ML2MC algorithm over the DLMC algorithm by combining the arguments in Theorem 3.11 of \cite{LMR1}, with the speedup estimate \eqref{eq:sp1}. 
\section{Extensions and Implementation}
\label{sec:5}
The multi-level machine learning algorithm \ref{alg:ML} can be readily extended in the following directions,
\subsection{Quasi-Monte Carlo type algorithms}
In the recent paper \cite{LMR1}, the authors obtained significantly lower generalization errors by choosing \emph{low-discrepancy sequences}, instead of randomly distributed points,  as the training set $\train$ in Step $1$ of the deep learning algorithm \ref{alg:DL}. The intuitive reason for this was the fact that low-discrepancy sequences, such as the Sobol and Halton points popularly used in Quasi-Monte Carlo integration methods \cite{CAF1}, are equi-distributed in the underlying parametric domain $Y$, see the recent paper \cite{MR1} for a rigorous explanation of this observation. We can readily adapt the multi-level algorithm \ref{alg:ML} to this setting by requiring that the training sets $\train_\ell$ for $0 \leq \ell \leq L$ are chosen as consecutive Sobol (or Halton) points. The rest of the algorithm is unchanged. Similarly, the UQ algorithm \ref{alg:ML2MC} can be readily adapted to this context.
\subsection{Other surrogate models.}
Deep neural networks are only one possible machine learning surrogate for the parameters to observable map $\map$ \eqref{eq:ptoob1}. Another popular class of surrogate models are \emph{Gaussian process regressions} \cite{GP}, which belong to a larger class of  so-called Bayesian models. Gaussian process regressions (GPR) rely on the assumption that the underlying map $\Delta$ is drawn from a Gaussian measure on a suitable function space, parameterized by,
\begin{equation}
\label{eq:gpr1}
\map({y}) \sim GP(m({y}), k({y}, {y}')),
\end{equation}
Here, $m({{y}})=\E{\map({y})}$ is the mean and $	k({{y}}, {y}')=\E{(\map({y}) - m({y})) (\map({y}') - m({y}'))}$ is the underlying covariance function. The mean and the covariance parametrize the so-called \emph{prior} measure. It is common to assume that $m \equiv 0$.

Given a training set $Y \supset \train = \{y_i\},1 \leq i \leq n$, the key idea underlying a Gaussian process regression is to apply Bayes theorem and update the conditional distribution for a test set $Y_{\ast} \subset Y, ~Y_{\ast} \cap \train = \emptyset$ with $\#(Y_{\ast}) = n_{\ast}$. For any $y_{\ast} \in Y_{\ast}$, denote $z_{\ast} = \map(y_{\ast})$, then one uses the Gaussian nature of the distributions to calculate the \emph{posterior} conditional probability by the formula,
\begin{equation}
    \label{eq:gpr2}
	{\rm Prob}({z}_*|\text{y}_*, \train) \sim \mathcal{N}\big(~G_*^T G^{-1}\text{z},~~ G_{**} - G_*^TG^{-1}G_*~\big),
\end{equation}
with $\mathcal{N}$ denoting a Gaussian distribution and $ z= [\map(y_1),\ldots,\map(y_n)]$, $y_i\in\train$. Here, $G\in\mathbb{R}^{n\times n}$, $G_{*}\in\mathbb{R}^{n\times n_*}$ and $G_{**}\in\mathbb{R}^{n_*\times n_*}$ are the training, the training-test and test Gram matrices, respectively, given by, 
\begin{equation}
\label{eq:gpr3}
G^{(i,j)} = k({y}_{i}, {y}_{j}),\quad G_{*}^{(i,j)} = k({y}_{i}, {y}_{*,j}),\quad{and}\quad G_{**}^{(i,j)} = k({y}_{*,i}, {y}_{*,j}).
\end{equation}
Thus, computation of the conditional probability \eqref{eq:gpr2} requires the inversion of the training Gram matrix $G$ \eqref{eq:gpr3}, for instance by a Cholesky algorithm, entailing a computational cost of ${\mathcal O}(n^3)$. 

Popular choices for the covariance function in \eqref{eq:gpr1} are the squared exponential (RBF function) and Matern covariance functions, 
\begin{equation}
\label{eq:SE_cov}
	k_{SE}({y}, {y}') = \exp{\bigg(-\frac{|| {y} - {y}'||^2}{2\ell^2}\bigg)}, \quad
	 k_{Matern}(y, y')={\frac {2^{1-\nu }}{\Gamma (\nu )}}{\Bigg (}{\sqrt {2\nu }}{\frac {||y -y'||}{\ell }}{\Bigg )}^{\nu }K_{\nu }{\Bigg (}{\sqrt {2\nu }}{\frac {||y -y'||}{\ell }}{\Bigg )}.
\end{equation}
Here $||\cdot||$ denotes the standard euclidean norm, $K_\nu$ is the Bessel function and $\ell$ the \textit{characteristic length}, describing the length scale of the correlations between the points $\text{y}$ and $\text{y}'$. 

We can readily adapt the multi-level algorithm \ref{alg:ML} to other surrogate models such as Gaussian process regressions. In fact, we propose a significantly more general form of the algorithm \ref{alg:ML} below,
\begin{algorithm} 
\label{alg:ext_ML} {\bf Multi-level learning of parameters to observable map}. 
\begin{itemize}
\item [{\bf Inputs}:] Parameterized PDE \eqref{eq:ppde}, Observable \eqref{eq:obsp}, high-resolution numerical method for solving \eqref{eq:ppde} and calculating \eqref{eq:obsp}, sequence of mesh resolutions $s_n$, number of training samples $N_0$, number of training samples $N_L$.
\item [{\bf Goal}:] Compute a machine learning surrogate $\map^{\ast}_{ml}$ for approximating the parameters to observable map $\map$ \eqref{eq:ptoob}. 

\item [{\bf Step $1$}:] For the coarsest mesh resolution $\ell_0$, select a training set $\train_0 = \{y^0_i\},~1 \leq i \leq N_0 = \#(\train_0)$, constituted by either random i.i.d points or with consecutive low-discrepency sequences such as Sobol points. For each $y^0_i$, compute $\map^{\Delta}_0(y^0_i)$ by solving the PDE \eqref{eq:ppde} on a mesh resolution $\Delta_0$ and computing the observable \eqref{eq:obsp}.  With the training data $\{y^0_i,\map^{\Delta_0}(y^0_i)\}_{1 \leq i \leq N_0}$, train the neural network $\map^{\ast}_{0,NN} \approx \map^{\Delta}_0$ by algorithm \ref{alg:DL} and the Gaussian Process regressor $\map^{\ast}_{0,GP} \approx \map^{\Delta}_0$, with suitable choice of the model hyperparameters. Assemble the ensemble model $\map^{\ast}_{0} = \alpha^0_{NN}\map^{\ast}_{0,NN} + \alpha^0_{GP}\map^{\ast}_{0,GP}$

\item [{\bf Step $2$}:] For each $\ell_k \in s_n$, $1 \leq k \leq n$, select a training set $\train_k = \{y^{\ell_k}_i\},~1 \leq i \leq N_k = \#(\train_k)$, with $N_k$ defined as in \eqref{eq:snumb} and $\train_k$ consisting of either random points or Sobol points. For each $y^{\ell_k}_i$, compute $\dtl_k(y^{\ell_k}_i)= \map^{\Delta_{\ell_k}}(y^{\ell_k}_i) - \map^{\Delta_{{\ell_{k-1}}}}(y^{\ell_{k}}_i)$ by solving the PDE \eqref{eq:ppde} on two successive mesh resolutions of $\Delta_{\ell_k}$ and $\Delta_{\ell_{k-1}}$ and computing the observable \eqref{eq:obsp}.  With the training data $\{y^{\ell_k}_i,\dtl_k(y^{\ell_k}_i)\}_{1 \leq i \leq N_k}$
train the neural network $\dtl^{\ast}_{k,NN} \approx \dtl_k$ and the Gaussian Process regressor $\dtl^{\ast}_{k,GP} \approx \dtl_k$, with suitable choice of the model hyperparameters. Assemble the ensemble model $\dtl^{\ast}_{k} = \alpha^k_{NN}\dtl^{\ast}_{k,NN} + \alpha^k_{GP}\dtl^{\ast}_{k,GP}$

\item [{\bf Step $3$}:] Form the machine learning surrogate $\map^{\ast}_{ml} \approx \map^{\Delta_L}$ as
\begin{equation}
    \label{eq:ml_ext}
    \map_{ml}^{\ast}(y) = \map^{\ast}_0(y) + \sum\limits_{k=1}^k \dtl^\ast_k(y), \quad \forall y \in Y.
\end{equation}

\begin{remark}
We have omitted for ease of notation the dependency of the trained network on the training set for ease of notation.
\end{remark}

\end{itemize}
\end{algorithm}
\subsection{Selection of hyperparameters.}
There are quite a few hyperparameters in the multi-level algorithm \ref{alg:ext_ML}. We choose these hyperparameters with the following procedure,
\subsubsection{Choice of Neural network hyperparameters}
The deep learning algorithm \ref{alg:DL} requires specification of the following hyperparameters: the architecture of the neural network \eqref{eq:ann1} (number of layers (depth) and size of each layer (width), the exponent $p$ in the loss function \eqref{eq:lf1}, the exponent $q$ and constant $\lambda$ in the regularization term \eqref{eq:lf2}, the choice of the optimization algorithm for minimizing \eqref{eq:lf2} and the starting value for it. 

Following \cite{LMR1}, we consider either $p=1$ or $p=2$ in the loss function \eqref{eq:lf1}. Similarly, the ADAM version of stochastic gradient algorithm \cite{ADAM} is used in \emph{full batch} mode with a learning rate of $\eta=0.01$ and is terminated at $10000$ epochs. For the regularization terms and starting value for ADAM, we use the \emph{ensemble training procedure} of \cite{LMR1} with either $q=1$ or $q=2$ and $\lambda = 5\times10^{-7}, 10^{-6}, 5\times10^{-6}, 10^{-5}, 5\times10^{-5}, 10^{-4}$. Similarly, five starting values (based on the He initialization) are used. Once the ensemble is trained, we select the hyperparameters that correspond to the smallest validation error, calculated on a validation set, created by setting aside $10 \%$ of the training samples. 

\subsubsection{Choice of Gaussian process hyperparameters.} For the characteristic length $l$ in \eqref{eq:gpr3}, we minimize the log negative marginal likelihood \cite{GP}. This leaves two hyperparameters i.e, the choice of the covariance kernel in \eqref{eq:gpr2} and if we choose the Matern kernel \eqref{eq:gpr3}, the choice of the parameter $\nu$. Here, we consider four hyperparameters configurations for the Gaussian process regression, namely either the squared exponential kernel in \eqref{eq:gpr3} or the Matern kernel with $\nu = 0.5,1.5,2.5$. These parameters are determined from an ensemble training process by selecting the parameter with the least validation error. 

\subsubsection{Choice of ensemble model coefficients.}
If $N_k>500$, $0\leq k\leq n$, we determine the ensemble hyperparameters $\alpha^k_{NN,GP}$ in algorithm \ref{alg:ext_ML}  by performing a linear least squares regression on the ensemble model with respect to a validation set, accounting for 10\% of the training set. Otherwise, we simply set $\alpha^k_{GP} = \alpha^k_{NN} = 0.5$.
\subsubsection{Choice of multi-level hyperparameters} 
In addition to the above hyperparameters, the multi-level algorithms \ref{alg:ML} and \ref{alg:ext_ML} involve three additional hyperparameters i.e, the sequence $s_n = \{\ell_k\}_{k=1}^n$ of mesh resolutions, parametrized by the model complexity $c_{ml}$ \eqref{eq:cml}, the number of samples $N_0$ at the coarsest level to learn $\map^{\Delta_0}$ and the number of samples $N_L$ to learn the detail $\dtl_n = \map^{\Delta_{\ell_n}} - \map^{\Delta_{\ell_{n-1}}}$, at the finest resolution. We will perform a sensitivity study to assess influence of these hyperparameters on the quality of results. 
 \begin{figure}[htbp]
\centering
\includegraphics[width=8cm]{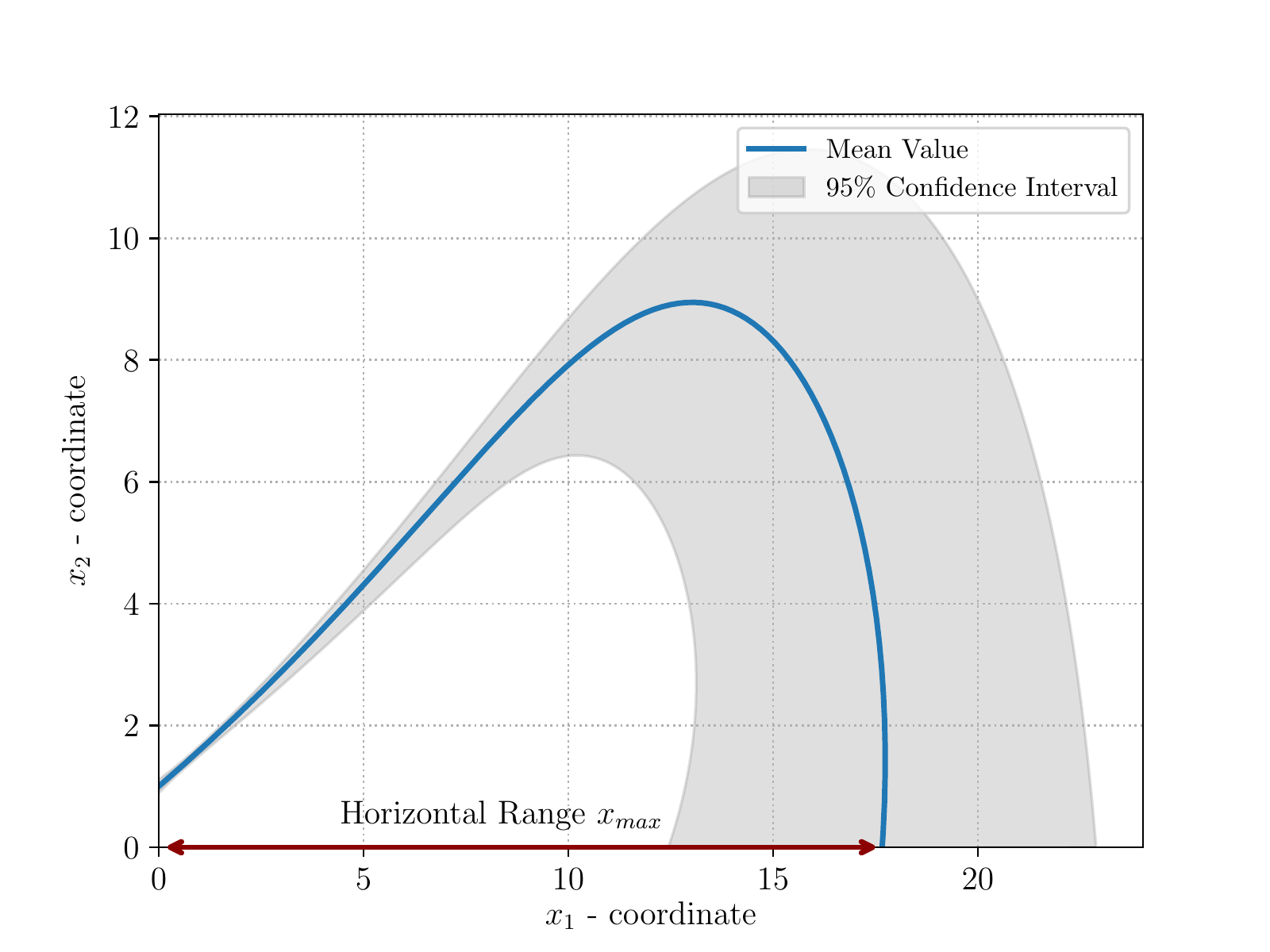}
\caption{An illustration of two-dimensional projectile motion, with the mean value and the envelope of trajectories, corresponding to the $95\%$ interval shown.}
\label{fig:pm1}
\end{figure}

\section{Numerical Experiments.}
\label{sec:6}
\subsection{Projectile motion.}
\label{sec:6.1}
We start with a dynamical system modeling the motion of a projectile, subjected to both gravity as well as air drag and described by the nonlinear system of ODEs,
\begin{equation}
\label{parab_syst}
\begin{aligned}
\frac{d}{dt}\boldsymbol{x}(t;y) &= \boldsymbol{v}(t;y), \quad 
\frac{d}{dt}\boldsymbol{v}(t; y) = - {{F}_D\big(\boldsymbol{v}(t;y); y \big)}\boldsymbol{e}_1 - g\boldsymbol{e}_2
\\
\boldsymbol{x}(y;0)&= \boldsymbol{x}_0(y), \quad 
\frac{d\boldsymbol{x}(y;0)}{dt} = \boldsymbol{v}_0(y).
\end{aligned}
\end{equation}
Here, ${F}_D = \frac{1}{2m}\rho C_d {\pi r^2}|| \boldsymbol{v}||^2$ denotes the drag force, with $\rho$ being the air density and $m$, $C_d$, $r$, the mass, the drag coefficient and the radius of the object, respectively.
Let further $\boldsymbol{x}_0(y) = [0, ~h]$, $\boldsymbol{v}_0(y) = [v_0\cos(\alpha), ~v_0\sin(\alpha)]$ be the initial position and velocity of the object (see figure \ref{fig:pm1} for a schematic representation). 

On account of measurement errors, the system is described by the following uncertain parameters,
\begin{equation}
\label{eq:pm2}
\begin{aligned}
\rho(y)&=1.225\big(1 + \varepsilon G_1(y)\big),~ 
r(y)=0.23\big(1 + \varepsilon G_2(y)\big),~ 
C_D(y)=0.1\big(1 + \varepsilon G_3(y)\big)
\\
m(y)&=0.145\big(1 + \varepsilon G_4(y)\big),~
h(y)=\big(1 + \varepsilon G_5(y)\big),~
\alpha(y)=30^\circ\big(1 + \varepsilon G_6(y)\big), ~
v_0(y)=25\big(1 + \varepsilon G_7(y)\big).
\end{aligned}
\end{equation}
Here, $y \in [0,1]^7$ describes the input parameter space, with uniform distribution, and $G_k(y)= 2y_k -1$ for $k=1,...,7$, with $\varepsilon = 0.1$. The objective of the simulation is to compute and quantify uncertainty with respect to the observable corresponding to the \textit{horizontal range} $x_{max}$ (see figure \ref{fig:pm1}), 
\begin{equation}
\label{eq:pm3}
    \map(y) = x_{max}(y) = x_1(y;t_f),\quad \text{with}~~ t_f = x_2^{-1}(0).
\end{equation}
Although at first glance, this problem appears simplistic, it has all the features of much more complicated problems (such as the one considered in the next experiment). Namely, the input parameter space $Y$ is moderately high-dimensional with $7$ dimensions and the parameters to the observable map is highly non-linear. Moreover, it yields a significant amount of variance in the trajectories and the observable (see figure \ref{fig:pm1}). Therefore, we will approximate solutions with a forward Euler discretization of \eqref{parab_syst}. The main advantage for choosing this model lies in the fact that we can compute solutions of \eqref{parab_syst} for a very large number of samples (realizations of $y$ in \eqref{eq:pm2}) with minimal computational work. This allows us to compute reference solutions and test our algorithms carefully. 
\begin{figure}[htbp]
    \begin{subfigure}{.49\textwidth}
        \centering
        \includegraphics[width=1\linewidth]{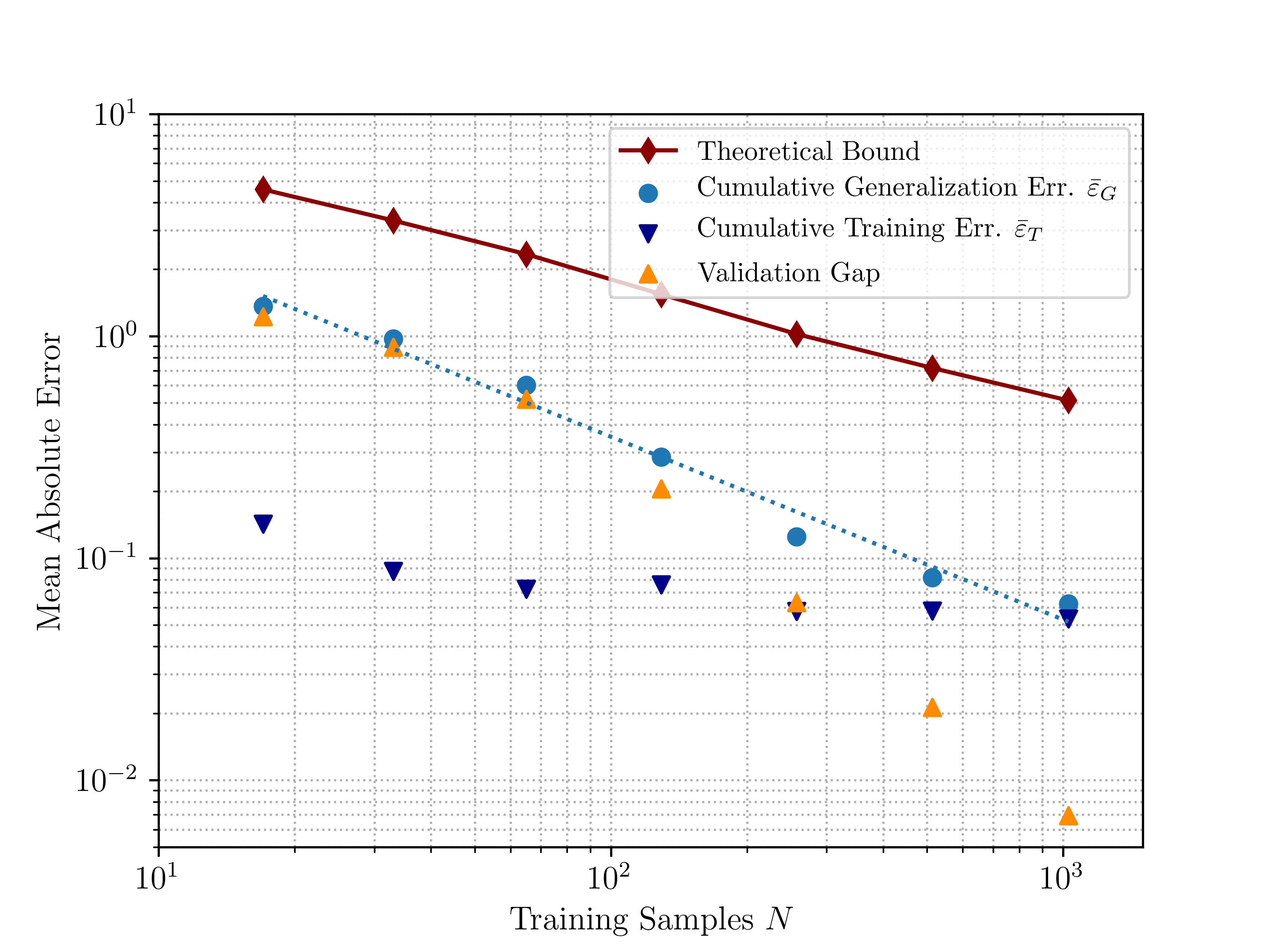}
        \caption{Errors.}
    \end{subfigure}
    \begin{subfigure}{.49\textwidth}
        \centering\
        \includegraphics[width=1\linewidth]{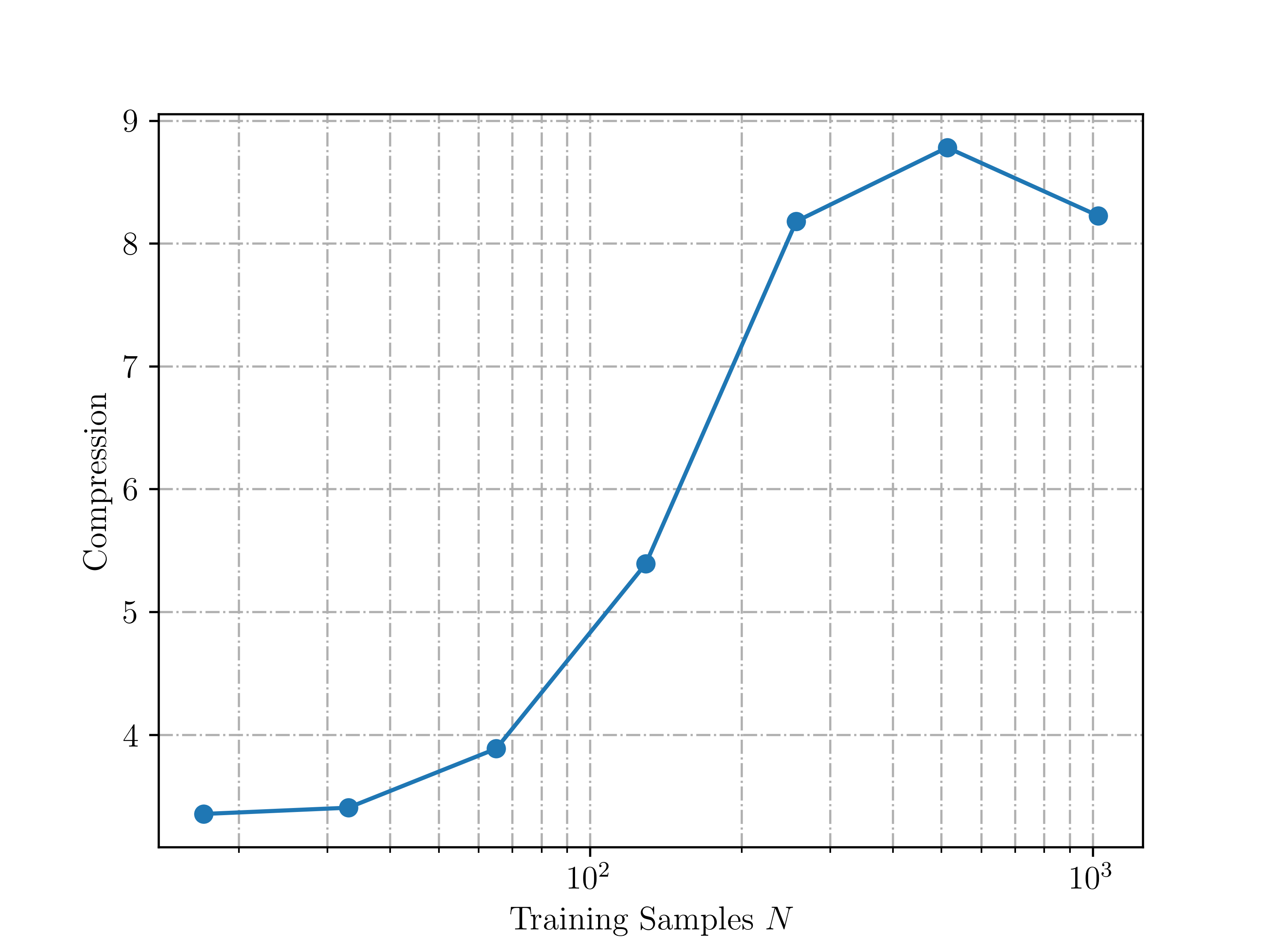}
        \caption{Compressions}
    \end{subfigure}
    \caption{Errors for the projectile motion. Left: Number of samples (X-axis) vs. training error, generalization error and computable upper bound \eqref{eq:genE2} on the generalization error. Right: Compression (Ratio of the upper bound to the generalization error) (Y-axis) vs. Number of training samples (X-axis).  }
    \label{fig:pm2}
\end{figure}

\subsubsection{Upper bounds on the generalization error.}
In section \ref{sec:23}, theorem \ref{theo:1}, we provided an upper bound \eqref{eq:genE}, \eqref{eq:genE2}, on the cumulative generalization error for the deep learning algorithm \ref{alg:DL}. Our subsequent theory about the utility of multi-level training, rested partially on this upper bound. We test our algorithm to check how sharp this upper bound is for this particular problem. 

To this end, we fix a neural network with $6$ hidden layers and $10$ neurons per layer, resulting in a network with $638$ tunable parameters (weights and biases). A mean absolute error loss function with a mean square regularization (with $\lambda = 10^{-6}$ in \eqref{eq:lf2}) is minimized with the ADAM optimizer for a fixed learning rate of $\eta =0.01$. To test the upper bound \eqref{eq:genE2}, we choose training sets as independent, uniformly distributed points in $[0,1]^7$ and the observable $\map$ \eqref{eq:pm3} is computed with a forward Euler discretization with time step $\Delta t = 0.00125$. 

We trained the above neural network with training sets of size $N_r=2^r,4 \leq r \leq 10$ and computed the generalization error with respect to $2000 - N_r$ i.i.d points in $[0,1]^7$. In order to compute the cumulative training, validation and generalization error, we retrain the network $K=60$ times and compute the corresponding average errors. The computation of the validation error is performed on a validation set of $N_r$ samples that is kept fixed over the $K$ resamplings of the training set. The expected validation gap \eqref{eq:vgap} is then computed by repeating the procedure above over $L=30$ different realizations of the validation set. 

As the right-hand side of \eqref{eq:genE2} involves computing standard deviations, we estimate ${\rm std}(\map)$ from $2000$ samples of $\map$ and ${\rm std}(\map^{\ast})$ from $1000$ realizations of the neural network. The same averaging strategy is also used for the computation of the variances.

The results are presented in figure \ref{fig:pm2}. In figure \ref{fig:pm2} (left), we plot the training error, validation gap, estimated generalization error and computed upper bound (rhs of \eqref{eq:genE2} versus the number of training samples. We see from this figure that the training error is consistently low and only decreases slightly as the number of training samples are increased. On the other hand, the validation gap is significantly larger than the training error for small number of samples but it decays very fast and is almost negligible when the number of training samples is large. On the other hand, the generalization error is approximately the same size as the validation gap (for small number of training samples) and the training error (for larger number of training samples). These findings are consisted with the theory presented here and infact, the generalization error decays as $N^{-0.82}$. On the other hand, the upper bound \eqref{eq:genE2} provides a reasonable overestimate of this generalization error. To further quantify the overestimate, we plot the \emph{compression} i.e, the ratio of the upper bound \eqref{eq:genE2} to the computed generalization error in figure \ref{fig:pm2} (right) and find that it lies between a factor of $3$ and $9$. This sharpness of the upper bound is particularly impressive given how difficult it is to obtain sharp upper bounds for neural networks \cite{ARORA,EMW1} and references therein. 

For the remaining part of this work we will consider only one single realization of the training set for ease of computation.
\begin{figure}[htbp]
    \centering
    \begin{subfigure}{.32\linewidth}
        \includegraphics[width=1\textwidth]{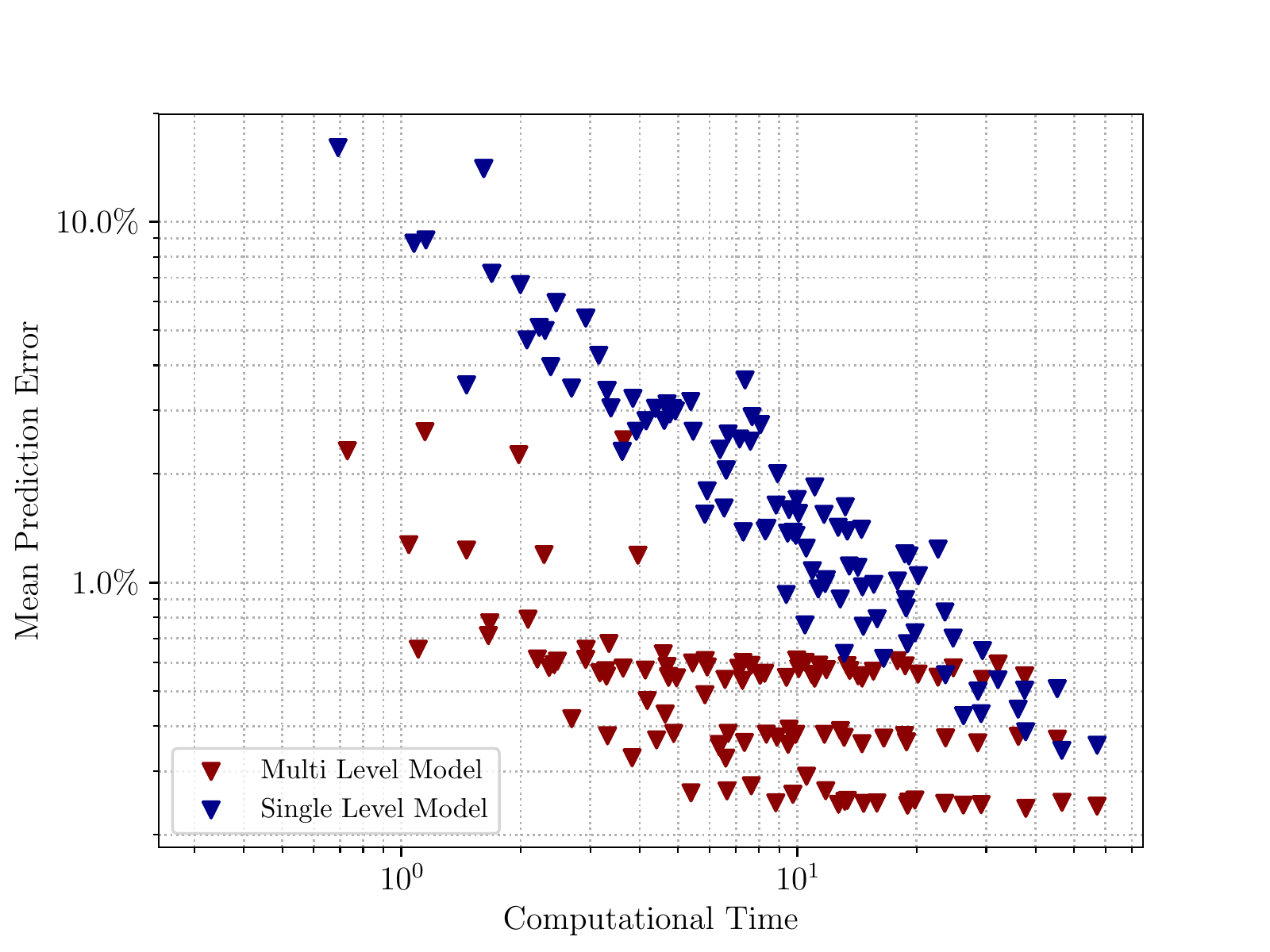}
        \caption{Prediction error vs. Computational cost}
    \end{subfigure}
    \begin{subfigure}{.32\linewidth}
        \includegraphics[width=1\textwidth]{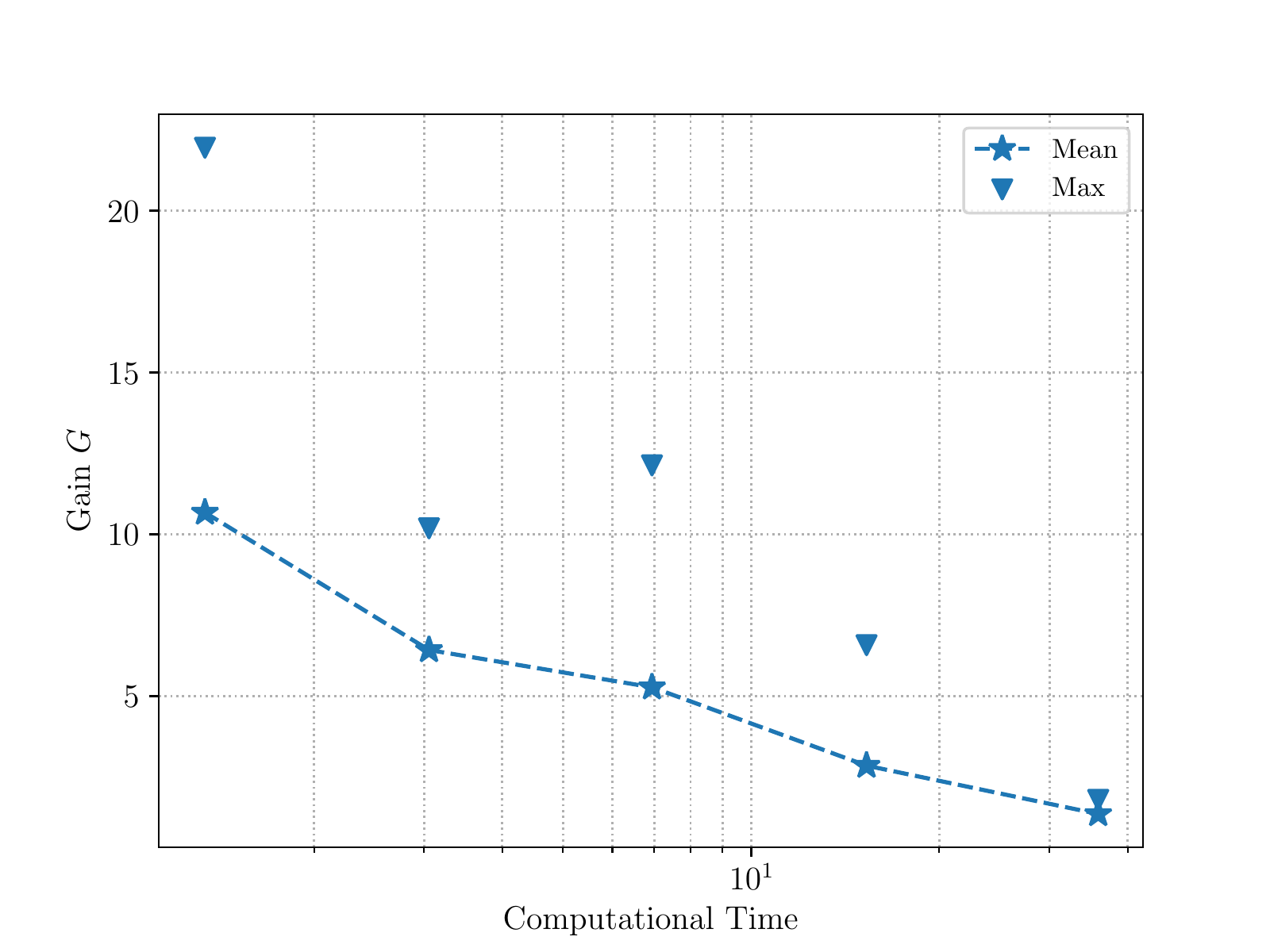}
        \caption{Gains \eqref{eq:gain1}}
    \end{subfigure}
    \begin{subfigure}{.32\textwidth}
        \includegraphics[width=1\linewidth]{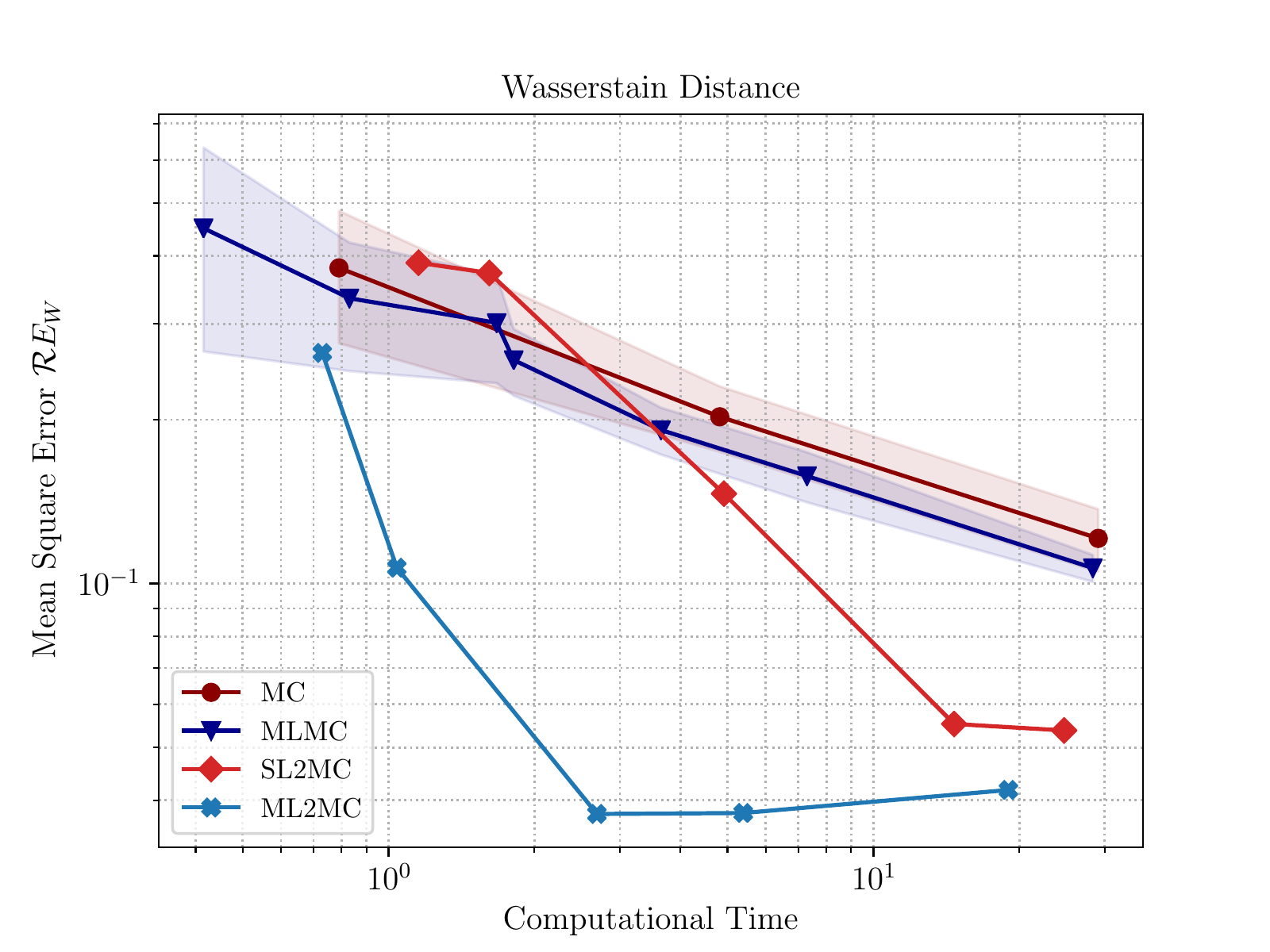}
        \caption{Wasserstein distance vs. Cost}
        \label{fig:x_max_wass}
    \end{subfigure}
    \caption{Results for the projectile motion case. Left: Prediction errors (Y-axis) vs computational cost (X-axis) for different multi-level parameters in algorithm \ref{alg:ML} and corresponding results with deep learning algorithm \ref{alg:DL}. Center: Gain \eqref{eq:gain1} vs. cost. Right: Errors, measured in the Wasserstein distance for the Monte-Carlo (MC), multi-level Monte Carlo (MLMC), single level machine learning (SL2MC) and multi-level machine learning algorithms (ML2MC) vs. Computational cost. }
    \label{fig:pm3}
\end{figure}
\subsubsection{Results of multi-level training and uncertainty quantification.}
Given the sharpness of the bound \eqref{eq:genE2} and the fact that the cumulative training error and validation gap seem to be significantly less than the generalization error (for sufficiently large number of training samples), we can follow Lemma \ref{lem:2} and expect that the multi-level machine learning algorithm \ref{alg:ML} will provide a lower generalization error than the deep learning algorithm \ref{alg:DL} for this problem.  To test this, we consider $7$ successive time step resolutions ($L=6$) within algorithm \ref{alg:ML}. Moreover, $\Delta t_0 = 0.08$ at coarsest resolution and $\Delta t_L = 0.00125$ at the finest resolution. We performed an ensemble training for the selection of the exponent $q$ of the regularization term and the parameter $\lambda$ of  $\map^\ast_{0}$ and $\dtl^{\ast}_{ref}$ (which the detail corresponding to the two coarsest resolutions) . The remaining parameters are identical to the previous subsection, with exception of the loss function, that in this case is the mean squared error ($p=2$ in \eqref{eq:lf1}).  $L^2$ regularization with $\lambda = 5\times 10^{-7}$ revealed the best performing configuration for both the maps. The same set of parameters found for $\dtl^{\ast}_{ref}$ was used for all the details. Moreover, we consider the following set of multi-level hyperparameters,
\begin{itemize}
    \item $4$ sequences of multi-level resolutions corresponding to $s^1_n = \{0, ~6\}$ ($c_{ml} = 0.16$), 
$s^2_n =\{0, ~3, ~6\}$ ($c_{ml} = 0.67$), 
$s^3_n = \{0, ~2, ~4, ~6\}$,~ $c_{ml} = 1.5$, $s^4_n = \{0, ~1, ~2, ~3, ~4,~5, ~6\}$ ($c_{ml} =6.0$).  
\item $4$ different choices of number of training samples at the coarsest resolution $N_0 = \{256, ~512,~ 1024, ~2048\}$
\item $7$ different choices of number of training samples at the highest resolution $N_L = \{4, ~8, ~16, ~32, ~64, ~92, ~128\}$
\end{itemize}

Overall, we have $112$ multi-level configurations. For each configuration, we run the multi-level algorithm \ref{alg:ML} and compute the generalization error from only high-fidelity samples i.e, samples at the finest mesh resolution. As we constrain the total number of high-fidelity samples to $2000$, we approximate the generalization error \eqref{eq:egen} with the prediction error ${\er}$,
\begin{equation}
\label{spe}
{\er}  := \frac{1}{N_{\test}}\Bigg(\sum_{y \in \test}\big |\map^{\Delta_L}(y) - {\map}^\ast(y)\big|^p \Bigg)^\frac{1}{p},
\end{equation} 
by choosing test sets  $\test \subset Y$ with number of samples $N_\test = \#(\test)$ ranging from $1872$ to $1996$. We choose $p=2$ in this section. For the sake of comparison, the deep learning algorithm \ref{alg:DL} is run on this data with the number of training samples determined by the need to match the cost of the multi-level algorithm \ref{alg:ML} (approximately). 

The corresponding results for the generalization error are plotted in figure \ref{fig:pm3} (Left). In this figure, we plot the generalization error for each multi-level hyperparameter configuration (and the corresponding equally expensive deep learning algorithm) versus the total computational time (in seconds). From this figure, we observe that the error with deep learning algorithm reduces with time (number of training samples). Moreover, the multi-level algorithms lead to a significant reduction in error at the same computational cost to the deep learning algorithm \ref{alg:DL}. A close inspection of this figure reveals that about $95\%$ of multi-level hyperparameters resulted in an error reduction (at the same cost) over the deep learning algorithm \ref{alg:DL} and about $7\%$ led to an order of magnitude reduction in error. 

This gain in efficiency is quantified in terms of the so-called \emph{gain}: \begin{equation}
    \label{eq:gain1}
    G = \frac{\er_{dl}}{\er_{ml}},
    \end{equation}
with $\er_{dl,ml}$ corresponding to the \emph{prediction errors} of the  the deep learning \ref{alg:DL} and multi-level \ref{alg:ML} algorithms, respectively. 
We plot the mean and the maximum value of the gain $G$  \eqref{eq:gain1} with respect to multi-level hyperparameters for a range of computational costs, in figure \ref{fig:pm3} (Center). From this figure, we see a mean gain between $2$ and $8$ and a maximum gain of $12$ for the multi-level algorithm over the deep learning algorithm \ref{alg:DL}. Larger gains were obtained for lower computational costs (less number of training samples), which is the case of practical interest. 

Finally, we apply the multi-level algorithm \ref{alg:ML} in the context of uncertainty quantification, in the form of the ML2MC algorithm \ref{alg:ML2MC}. We compute approximations $\hat{\mu}^{\ast}_{ml2mc}$ to the full push-forward measure $\hat{\mu}^{\Delta}$ \eqref{eq:pf1}. A reference push-forward measure $\hat{\mu}^{ref}$ is computed with a very small time step of $\Delta t = 0.001$ and $20000$ Monte Carlo samples and we ascertain the quality of the ML2MC algorithm by computing \emph{Wasserstein distances} ${\mathcal W}_1 \left(\hat{\mu}^{\ast}_{ml2mc},\hat{\mu}^{ref}\right)$ with the function \textit{wasserstein\_distance} of the Python library \textit{scipy.stats} \cite{py_wass}, on the multi-level hyperparameter configurations shown in Table \ref{tab:parab_MLmlearn}. These hyperparameters approximately correspond to those that result in the highest gain in the prediction error (see figure \ref{fig:pm3} (A)). 

In order to compare the multi-level UQ algorithm \ref{alg:ML2MC} with existing algorithms, we select the following,
\begin{itemize}
    \item Standard Monte Carlo approximation of the push forward measure, at the finest resolution of $\Delta t = 0.00125$.
    \item A multi-level Monte Carlo algorithm for computing push-forward measures as proposed in \cite{LyePhD}. 
    \item The \emph{single-level} variant (SL2MC) of algorithm \ref{alg:ML2MC} where the multi-level algorithm \ref{alg:ML} is replaced by the deep learning algorithm \ref{alg:DL}, trained with 9,~15,~21,~64,~191 and 322 number of samples at the finest mesh resolution. This single-level algorithm coincides with the DLMC algorithm of \cite{LMR1}.
\end{itemize}
In figure \ref{fig:pm3} (Right), we plot the mean Wasserstein distance versus the computational time for all four competing algorithms. We observe from this figure that the ML2MC algorithm \ref{alg:ML2MC}, clearly outperforms the competing algorithms. On an average (over the computational costs considered), it provides a speedup (reduction in error at same computational cost) of a factor of $5$ over the Monte Carlo algorithm, $4$ over the MLMC method and $3$ over the single-level machine learning UQ algorithm, based on the deep learning algorithm \ref{alg:DL}.

\begin{table}[htbp] 
    \centering 
    \footnotesize{
        \begin{tabular}{c c |  c |  c|  c | c } 
            \toprule
            \bfseries  &\multicolumn{5}{c}{\bfseries ML2MC Configurations}\\
            \midrule
            \bfseries Samples $N_0$  & 256 & 256 & 2048 & 2048 & 2048\\
            \bfseries Samples $N_L$ & 4 & 8 & 8 & 32 & 64\\
            
            \bfseries Complexity $c_{ml}$ & 0.17 & 0.67 & 0.67 & 0.67 & 1.5 \\
            
            \bottomrule 
        \end{tabular}
    }
        \caption{Different cost configurations used to evaluate the performance of multilevel machine learning method  for the Projectile Motion example.}
        \label{tab:parab_MLmlearn}
\end{table}
\begin{figure}[htbp]
    \begin{subfigure}{.36\textwidth}
        \centering
        \includegraphics[width=1\linewidth]{{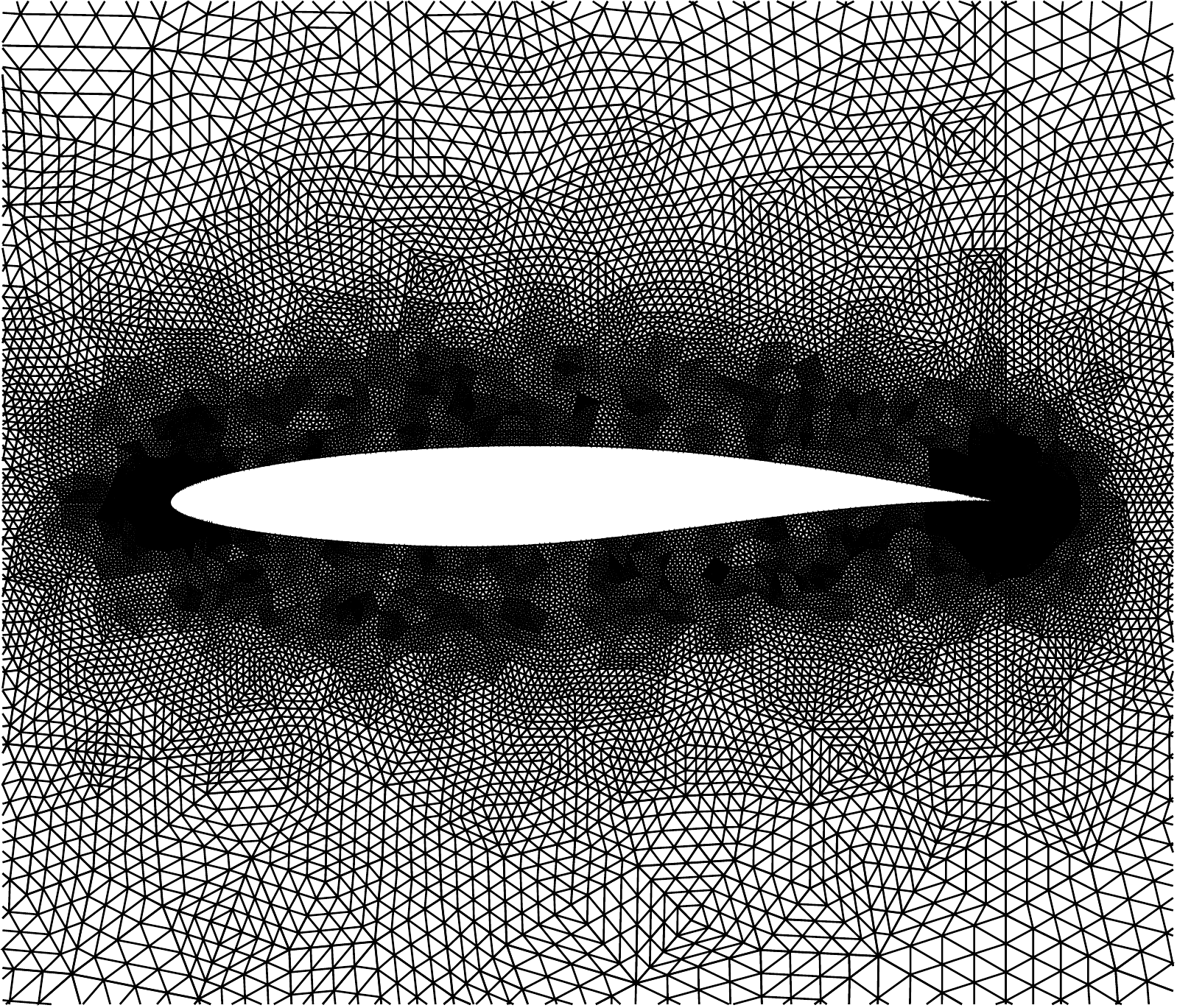}}
        \caption{Hi-resolution grid}
    \end{subfigure}
    \begin{subfigure}{.3\textwidth}
        \centering\
        \includegraphics[width=1\linewidth]{{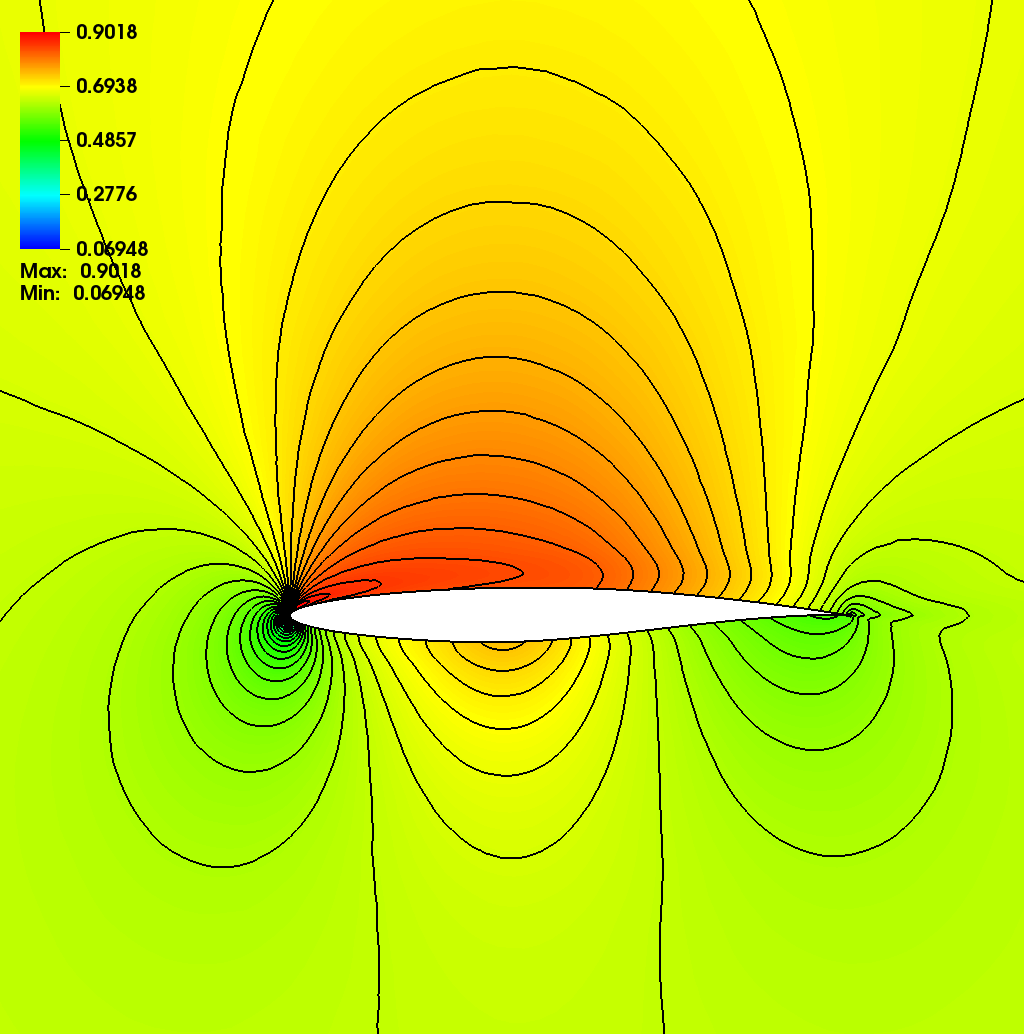}}
        \caption{Mach number (Sample)}
    \end{subfigure}
    \begin{subfigure}{.3\textwidth}
        \centering\
        \includegraphics[width=1\linewidth]{{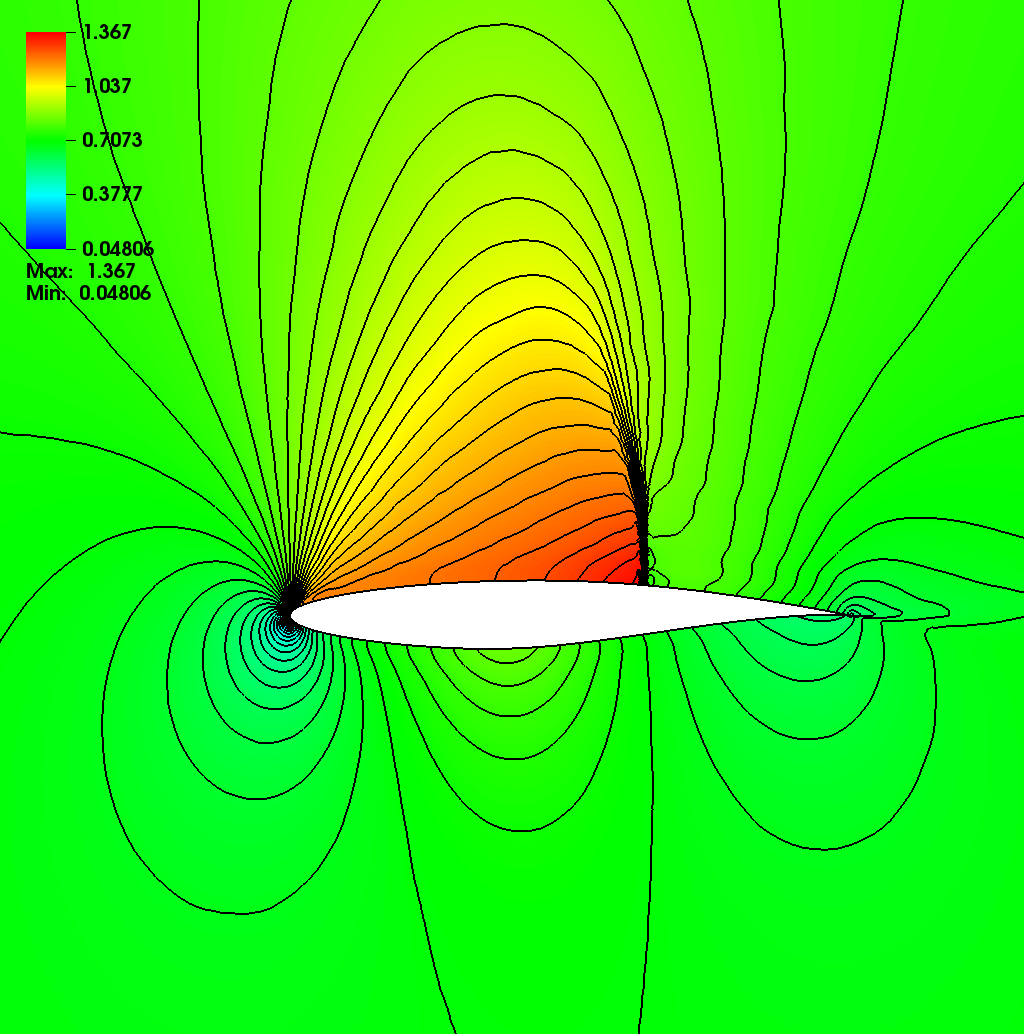}}
        \caption{Mach number (Sample)}
    \end{subfigure}
    
    \caption{Flow past a RAE2822 airfoil. Left: High-resolution grid Center and Right: Flow visualized with Mach number for two different samples. }
    \label{fig:fpa1}
\end{figure}

\subsection{Flows past airfoil.}
\label{sec:6.2}
In this section, we consider a much more realistic example of a compressible flow past a RAE2822 airfoil. The problem set up is a benchmark for UQ in fluid flows \cite{UMRIDA} and is identical to the one considered in \cite{LMR1}. The two-dimensional compressible Euler equations are solved on the following perturbed free stream conditions and the profile geometry:
\begin{equation}
\begin{aligned}
T^\infty(y)&= \big(1 + \varepsilon_1 G_1(y)\big),\quad 
M^\infty(y)=0.729\big(1 + \varepsilon_1 G_2(y)\big), 
\\[5pt]
p^\infty(y)&= \big(1 + \varepsilon_1 G_2(y)\big),\quad
\alpha(y)=2.31^\circ\big(1 + \varepsilon_1 G_6(y)\big),
\\[5pt]
S_L(x;y)&= \overline{S}_L(x)\big(1 + \varepsilon_2 G_4(y)\big),\quad 
S_U(x;y)=\overline{S}_U(x)\big(1 + \varepsilon_2 G_4(y)\big),
\end{aligned}
\end{equation}
where $\alpha$ is the angle of attack and $\bar{S}_U(x)$, $\bar{S}_L(x)$, $x\in[0,1]$, denote the unperturbed upper and lower surfaces of the airfoil, respectively, see figure \ref{fig:fpa1} (Left) for the reference geometry. $G_k$ are defined as in the previous numerical examples and $\varepsilon_1=0.1$, $\varepsilon_2=0.2$. 
The observables are the lift and drag coefficients:
\begin{equation}
\map_1(y) =
C_L(y) = \frac{1}{E_k^\infty(y)}\int_{S_L\cup S_U}p(y)n(y)\cdot\hat{y}(y)ds,\quad
\end{equation}
\begin{equation}
\map_2(y) =
C_D(y) = \frac{1}{E_k^\infty(y)}\int_{S_L\cup S_U}p(y)n(y)\cdot\hat{x}(y)ds,\quad
\end{equation}
where $\hat{x}(y) = [\cos(\alpha(y)), \sin(\alpha(y))]$, $\hat{y}(y) = [-\sin(\alpha(y)), \cos(\alpha(y))]$ and
\begin{equation}
    E_k^\infty(y)= \frac{\rho^\infty(y)||\textbf{u}^\infty(y)||^2}{2}
\end{equation}
is the free-stream kinetic energy.

Thus, the input parameter space is the $6$-dimensional cube $Y = [0,1]^6$ and two samples, corresponding to two different realizations in $Y$ are shown in figure \ref{fig:fpa1}. 

In \cite{LMR1}, the authors observed that for this problem, the deep learning algorithm \ref{alg:DL} with Sobol training points was significantly more accurate, with a factor of $10-20$ lower generalization errors than with randomly distributed training points. Hence, we will only consider the case of Sobol training points here. Moreover, we use this example to test the more general multi-level machine learning algorithm \ref{alg:ext_ML}. 

To this end, we consider $5$ levels of grid resolution approximating the flow past the airfoil, with $L=4$ in algorithm \ref{alg:ext_ML}. The grid at the finest resolution is shown in figure \ref{fig:fpa1} (Left). For each resolution, the two-dimensional Euler equations will be solved with the TEnSUM code, which implements vertex centered high-resolution finite volume schemes on unstructured triangular grids \cite{RCFM}. 

For implementing the multi-level algorithm, we consider $4$ sequences of multi-level resolutions corresponding to $s^1_n = \{0,4\}$ ($c_{ml} = 0.25$), 
$s^2_n =\{0,2,4\}$ ($c_{ml} = 1.0$), 
$s^3_n = \{0,2,3,4\}$ ($c_{ml} = 2.31$), $s^4_n = \{0,1,2,3,4\}$ ($c_{ml} =4.0$). Moreover, the same choices of samples at the coarsest level ($N_0$) and finest level ($N_L$) are made as in the previous numerical experiment. This results in a total of $112$ multi-level hyperparameter configurations. 
\begin{figure}[htbp]
    \begin{subfigure}{.45\textwidth}
        \centering
        \includegraphics[width=1\linewidth]{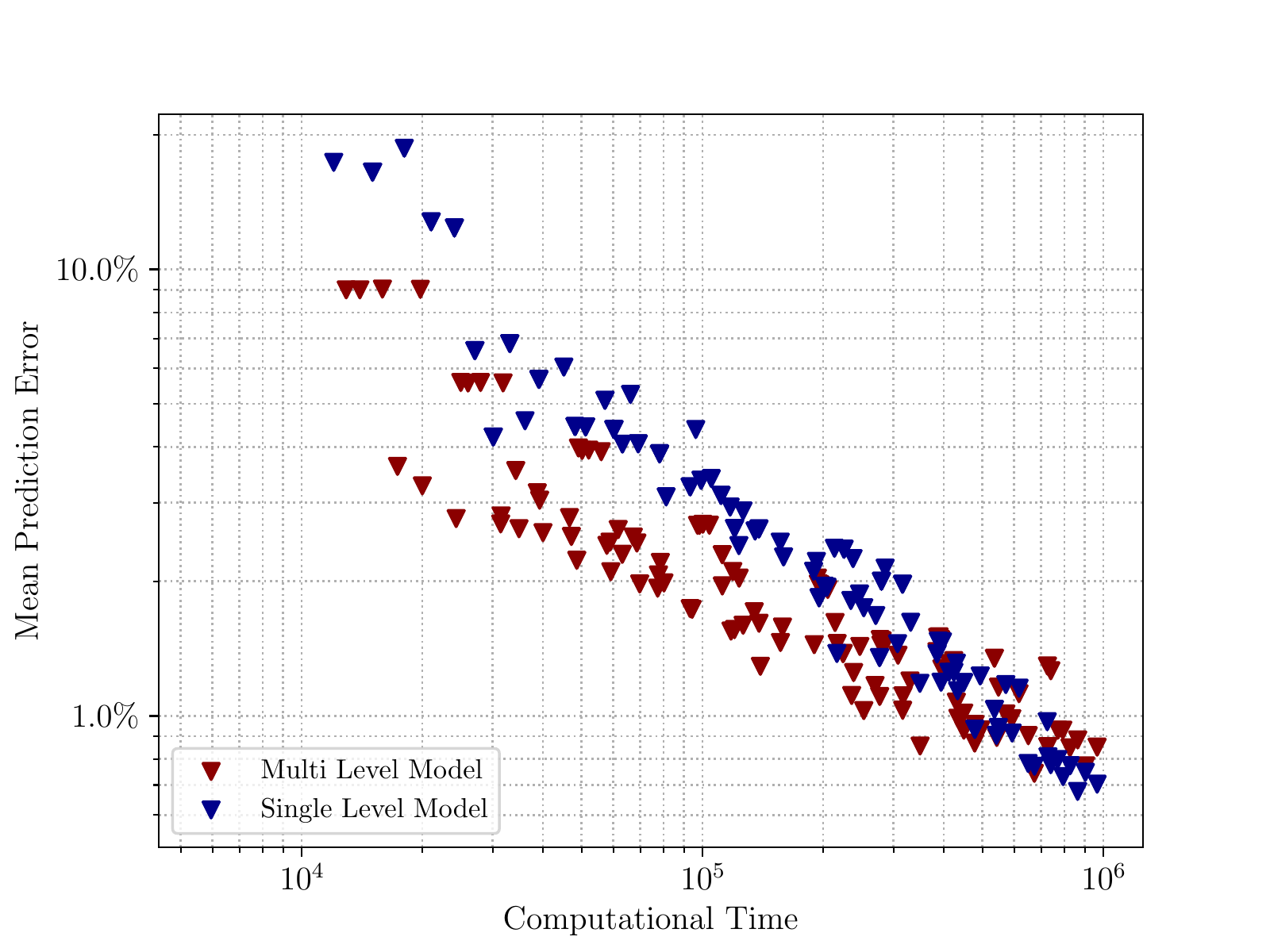}
    \end{subfigure}
    \begin{subfigure}{.45\textwidth}
        \centering\
        \includegraphics[width=1\linewidth]{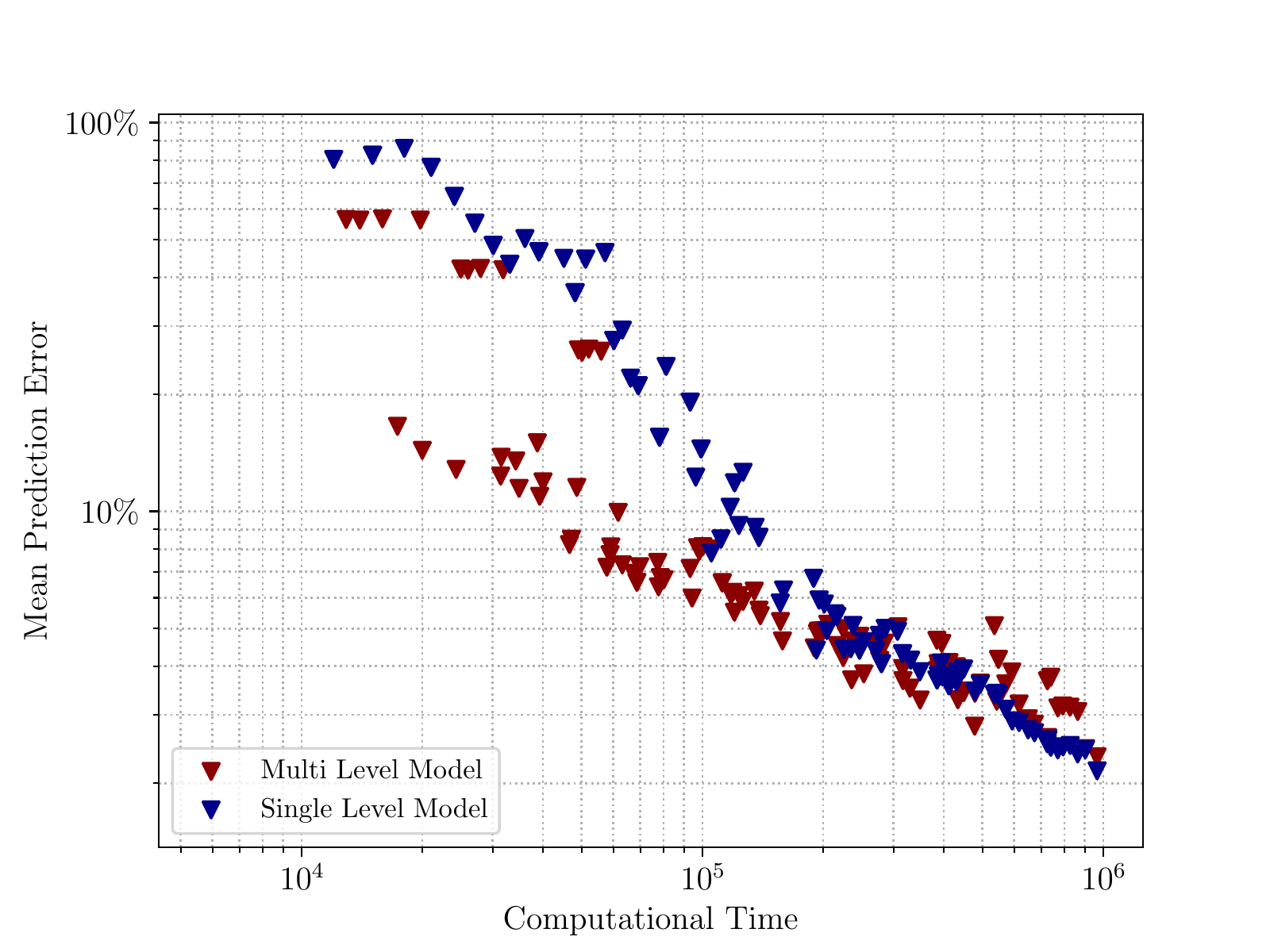}
    \end{subfigure}

    \begin{subfigure}{.45\textwidth}
        \centering\
        \includegraphics[width=1\linewidth]{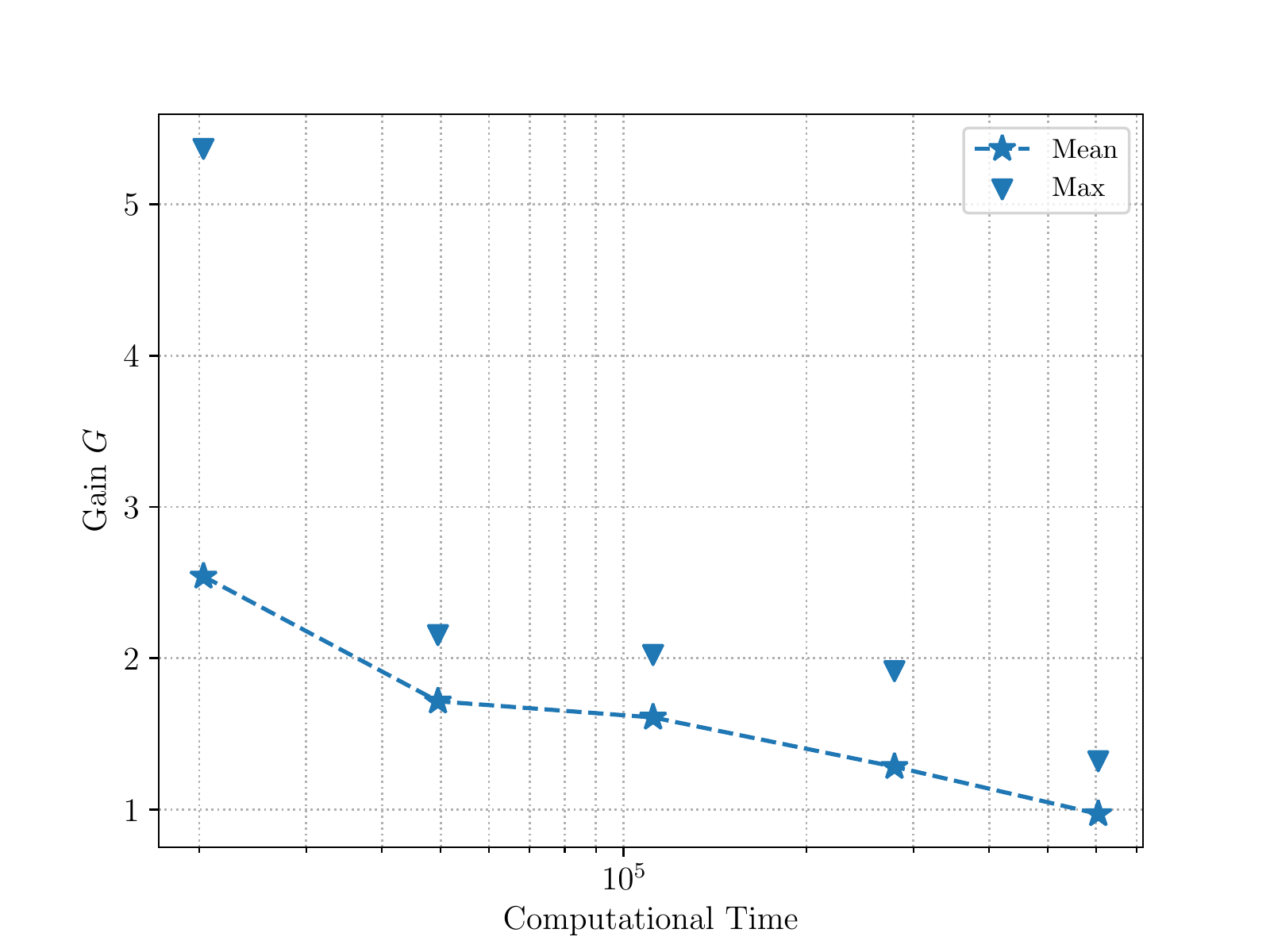}
        \caption{Lift}
    \end{subfigure}
        \begin{subfigure}{.45\textwidth}
        \centering\
        \includegraphics[width=1\linewidth]{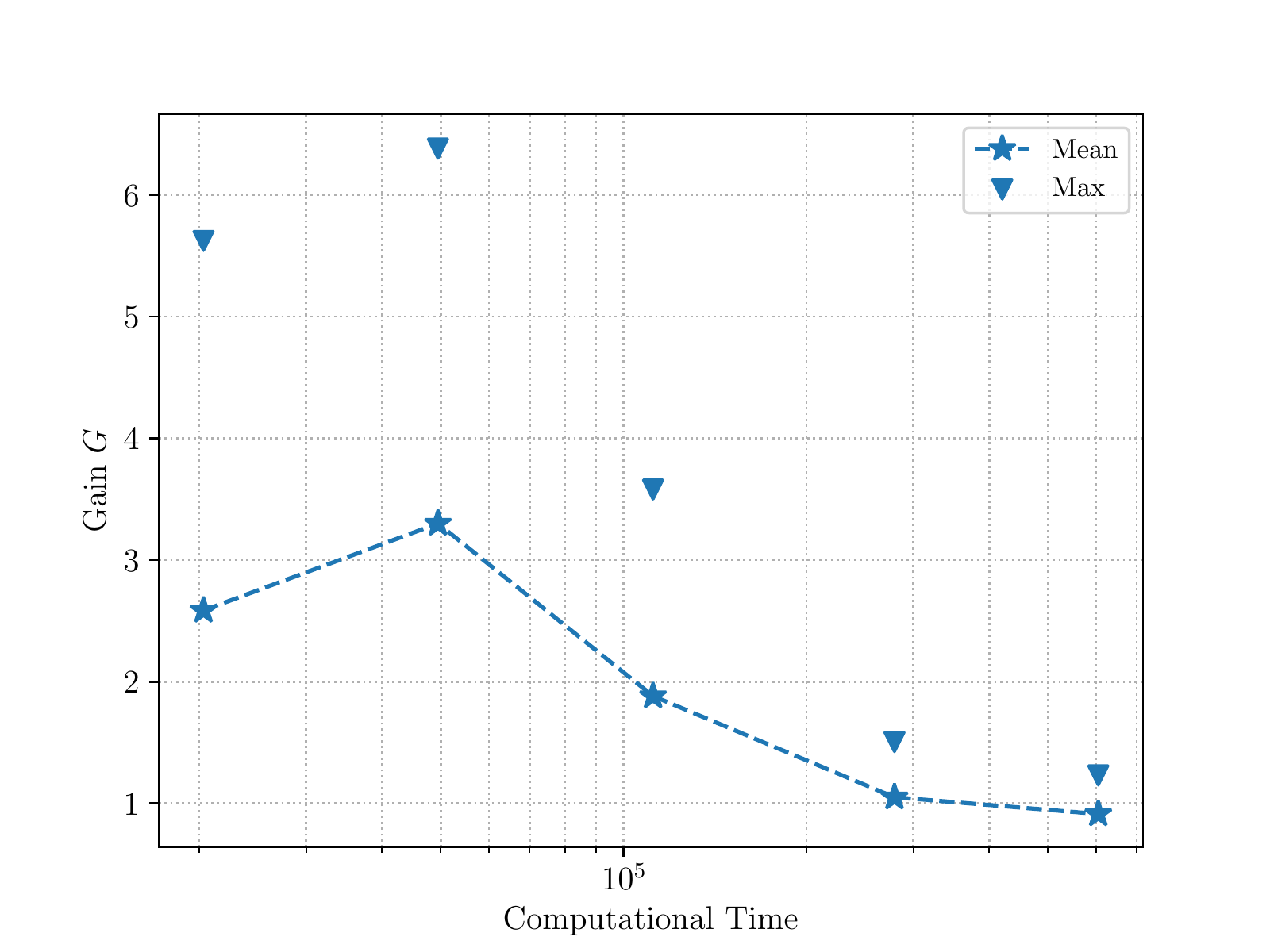}
        \caption{Drag}
    \end{subfigure}
    \caption{Prediction errors with the multi-level machine learning algorithm \ref{alg:ext_ML} for flow past a RAE2822 airfoil. Top Row: Prediction error vs. Computational cost. Bottom Row: Gain \eqref{eq:gain1} over single level algorithm vs. computational cost.}
    \label{fig:fpa2}
\end{figure}

For the lift coefficient, we use the neural network architecture of \cite{LMR1}, namely a fully-connected network of $9$ hidden layers with $12$ neurons in each layer, and minimize a mean square loss function with the ADAM optimizer at a fixed learning rate of $\eta = 0.01$. This configuration was used for learning both $\map^\ast_{0,NN}$ as well as the reference detail $\dtl^{\ast}_{ref,NN}$. An ensemble training, in the sense of \cite{LMR1}, was performed to discover that a mean-square regularization with a $\lambda = 5 \times 10^{-7}$, provided the optimal hyperparameters for both maps. The same hyperparameter configuration was used for learning other details. Similarly for the Gaussian process regression hyperparameters, we performed an ensemble training on $\map^{\ast}_{0,GP}, \dtl^{\ast}_{ref,GP}$ to find that the Matern covariance kernel \eqref{eq:gpr3} with $\nu = 1.5$ provided the best hyperparameters. 

For the drag coefficient, we retain the same hyperparameters for the Gaussian process, but for the neural networks in algorithm \ref{alg:ext_ML}, we performed a full ensemble training for each map, which resulted in hyperparameters shown in table \ref{setting_drag_net}. In addition to those listed in this table, the ADAM optimizer was used to minimize a mean-square loss function. \begin{table}[htbp] 
    \centering
    \renewcommand{\arraystretch}{1.1} 
    
    \footnotesize{
        \begin{tabular}{c  c c c c c c } 
            \toprule
            \bfseries   &\bfseries Samples & \bfseries Lear. Rate &\bfseries $L^1$ - Reg. & \bfseries $L^2$ - Reg. &\bfseries Depth &\bfseries Width \\ 
            \midrule
            \midrule
            $\map^{\Delta_0}$            & 1024   & 0.001 &0& $5\cdot10^{-7}$& 6& 16\\
            \midrule 
            $\map^{\Delta_1} - \map^{\Delta_0}$ & 256  & 0.001 &0& $10^{-6}$& 6& 12\\
                \midrule 
            $\map^{\Delta_2} - \map^{\Delta_1}$ & 64  & 0.001 &$5\cdot10^{-7}$& 0& 6& 12\\
                \midrule 
            $\map^{\Delta_3} - \map^{\Delta_2}$ & 64 & 0.01 &$5\cdot10^{-6}$& 0& 6& 8\\
                \midrule 
            $\map^{\Delta_4} - \map^{\Delta_3}$ & 8 &0.01 &$10^{-5}$& 0& 6& 8\\
                \midrule 
            $\map^{\Delta_2} - \map^{\Delta_0}$ & 256  & 0.001 &$5\cdot10^{-7}$& 0& 6& 8\\
            \midrule
            $\map^{\Delta_4} - \map^{\Delta_2}$ & 8  &0.01 &$5\cdot10^{-7}$& 0& 6& 12\\
                \midrule 
            $\map^{\Delta_3} - \map^{\Delta_0}$ & 256  & 0.001 &$5\cdot10^{-7}$& 0& 9& 16\\
            \midrule 
            $\map^{\Delta_4} - \map^{\Delta_0}$ & 32  &0.001 &$10^{-5}$& 0& 6& 8\\
            \bottomrule
        \end{tabular}
    \caption{Best performing \textit{Neural Network} configurations for the drag coefficient. In first column, we show the number of training samples used in performing the ensemble training.}
        \label{setting_drag_net}
    }
\end{table}
\subsubsection{Prediction errors with multi-level model.}
We approximate the generalization error \eqref{eq:genE} by computing the prediction error of the algorithms \eqref{spe} on a \emph{test set}, formed by $2000$ Sobol training points on the finest mesh resolution. The resulting prediction errors for the multi-level machine learning algorithm \ref{alg:ext_ML}, for all the $112$ multi-level hyperparameters configurations, for both the lift and the drag, are shown in figure \ref{fig:fpa2} (Top Row). To provide a comparison, we also compute the \emph{single-level} surrogate map $\alpha_1 \map^{\ast}_{NN} + \alpha_2 \map^{\ast}_{GP}$, with all the training samples being generated at the finest grid resolution and determined to (approximately) match the computational cost of the multi-level algorithm and plot the corresponding prediction errors in figure \ref{fig:fpa2}. From this figure, we see that the multi-level algorithms clearly outperform the single-level machine learning algorithm, with approximately $75 \%$ of the configurations, resulting in a lower prediction error at similar computational cost. A few of the configurations ($3 \%$ for the lift and $8 \%$ for the drag) result in at least a factor of $4$ reduction in the prediction error at same cost. 

This gain in efficiency is further quantified by computing the gain \eqref{eq:gain1} and plotting it in figure \ref{fig:fpa2} (Bottom Row). From this figure, we see that on an average, the gain with the multi-level model is a factor of $2-3$ for the Lift and $3-4$ for the Drag, in the range of reasonably small computational times, corresponding to the interesting case of low number of training samples. The maximum gain in this range is approximately $6$ for both observables. Although these gains are a bit smaller than the ones for the projectile motion, they are more impressive as the underlying maps are hard to learn but the combination of Sobol training points, deep neural networks, and Gaussian processes makes the competing single-level machine learning model quite accurate. 
\subsubsection{Sensitivity of results to multi-level hyperparameters.} 
\begin{figure}[htbp]
    \centering
    \begin{subfigure}{.325\linewidth}
        \includegraphics[width=1\textwidth]{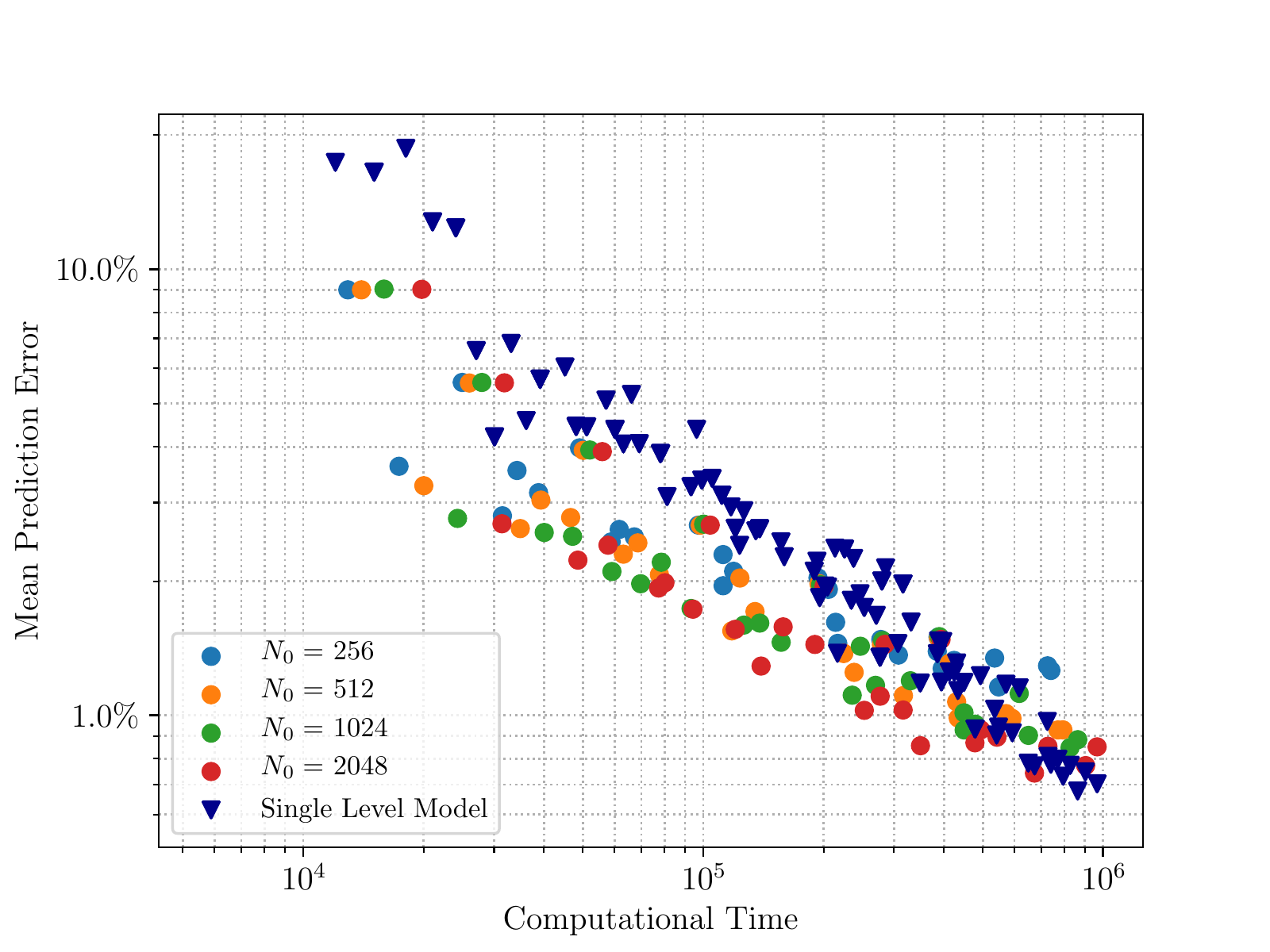}
    \end{subfigure}
    \begin{subfigure}{.325\linewidth}
        \includegraphics[width=1\textwidth]{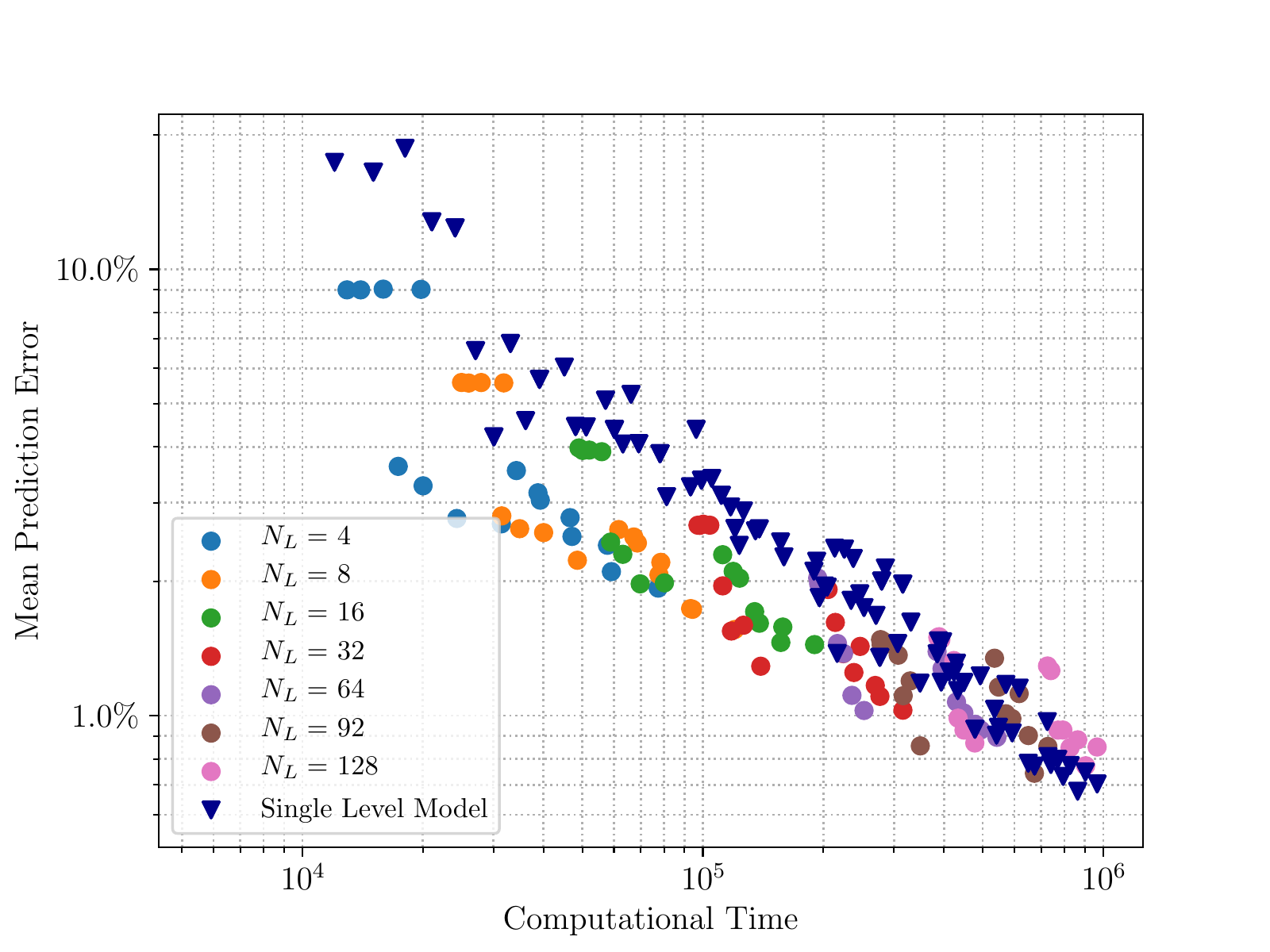}
    \end{subfigure}
    \begin{subfigure}{0.325\linewidth}
        \includegraphics[width=1\textwidth]{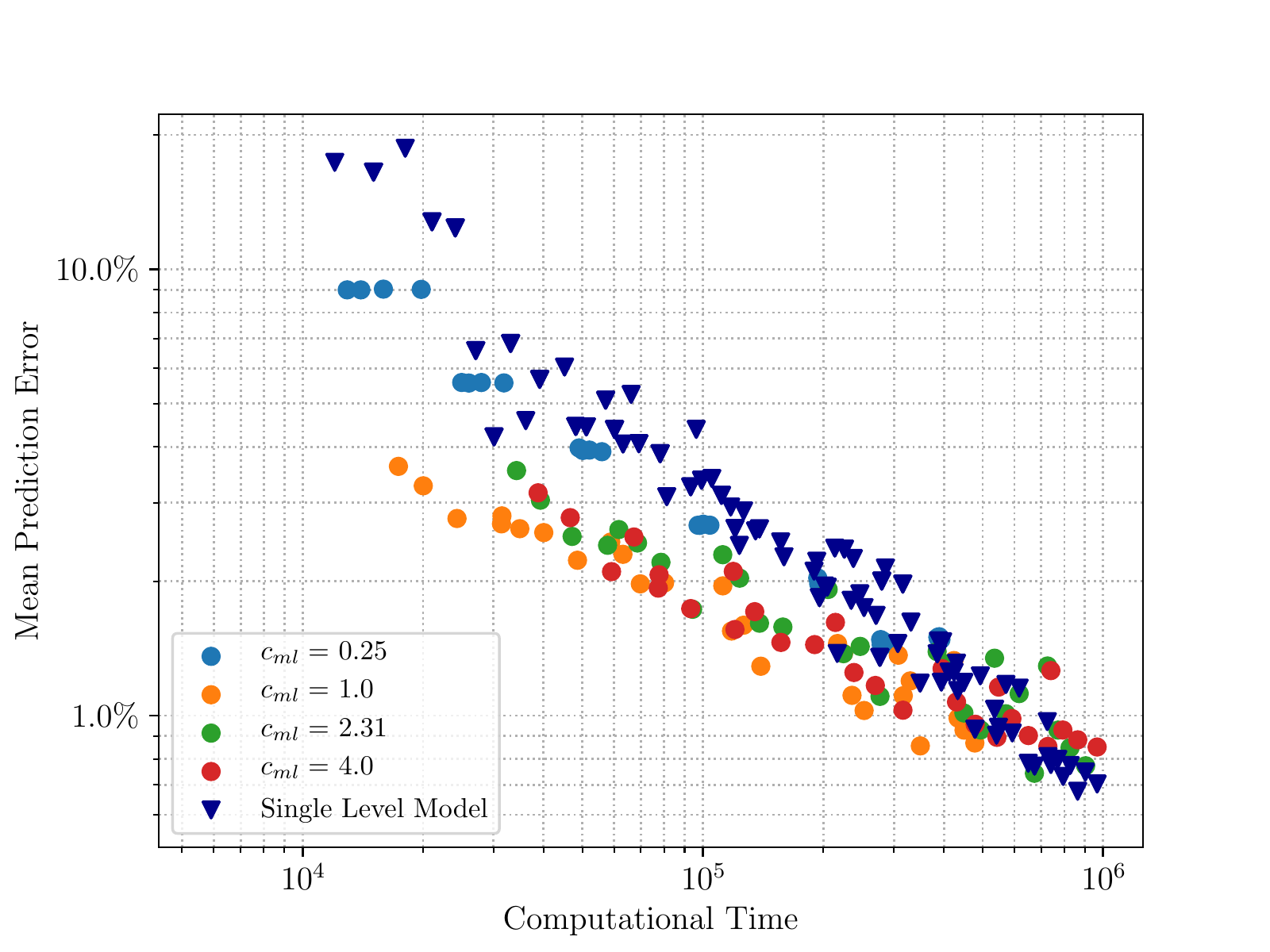}
    \end{subfigure}
    
    \begin{subfigure}{.325\linewidth}
        \includegraphics[width=1\textwidth]{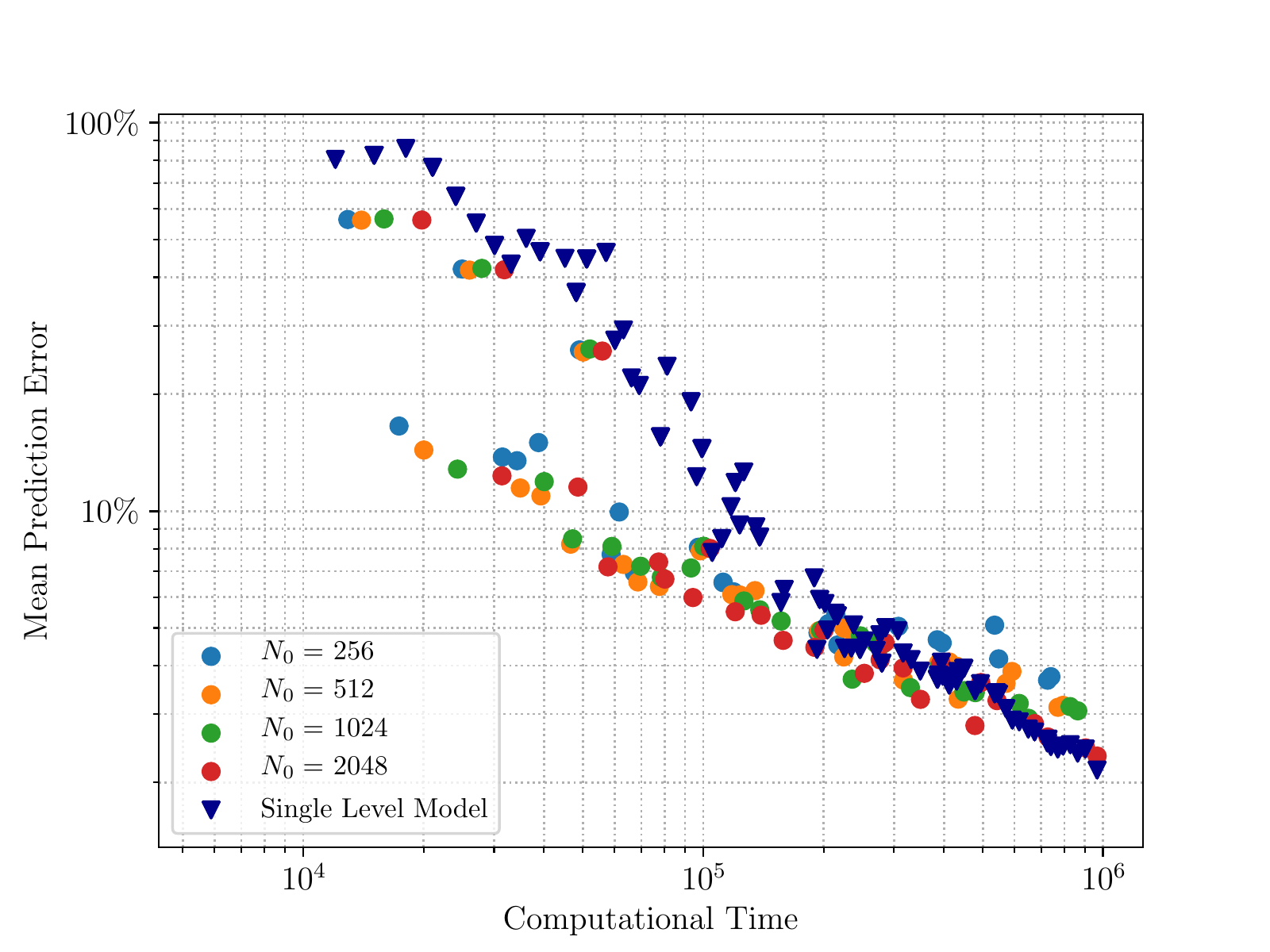}
        \caption{Sensitivity to $N_0$}
    \end{subfigure}
    \begin{subfigure}{.325\linewidth}
        \includegraphics[width=1\textwidth]{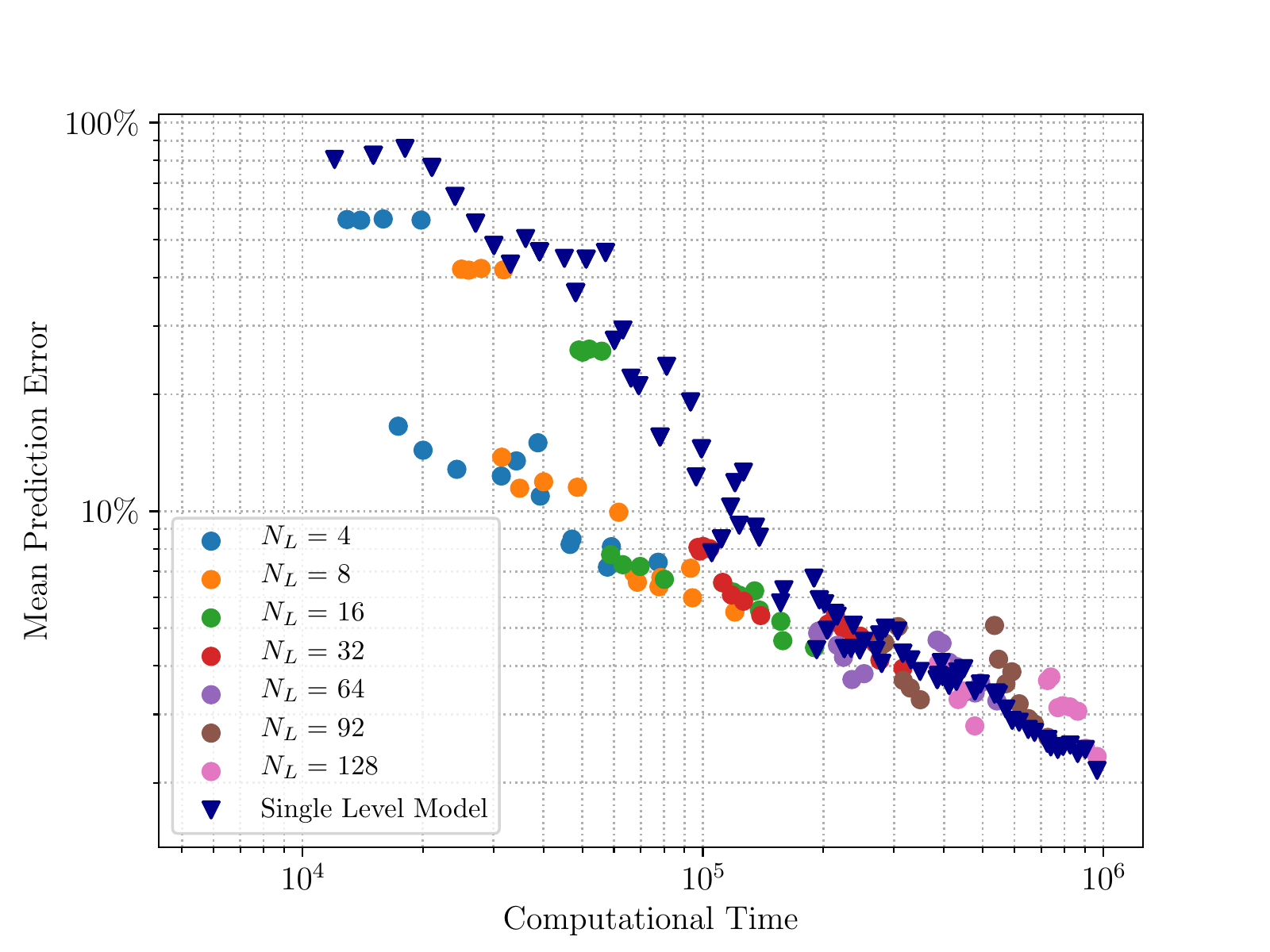}
        \caption{Sensitivity to $N_L$}
    \end{subfigure}
    \begin{subfigure}{0.325\linewidth}
        \includegraphics[width=1\textwidth]{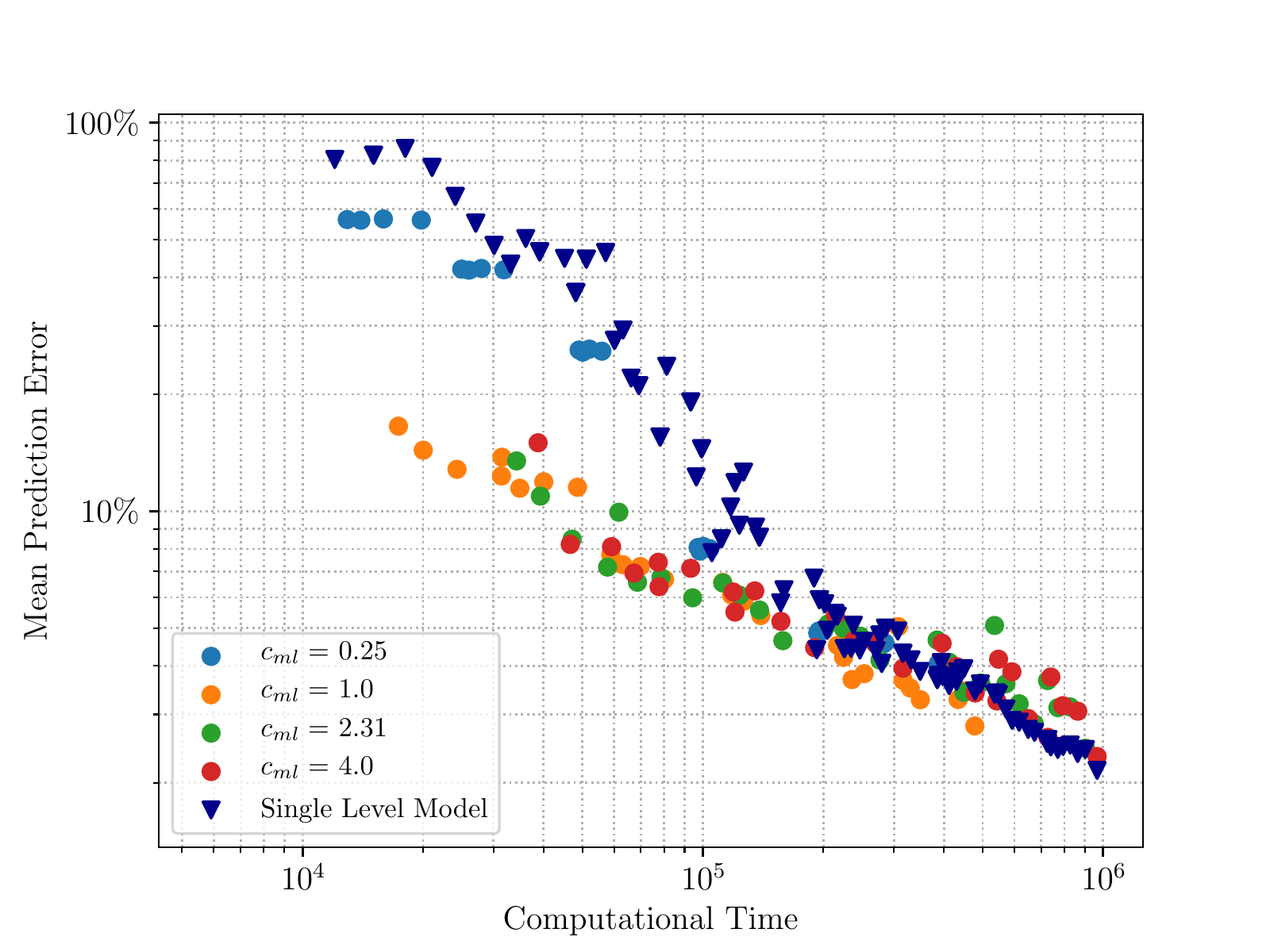}
        \caption{Sensitivity to $c_{ml}$}
    \end{subfigure}
\caption{Sensitivity of the multi-level machine learning algorithm \ref{alg:ext_ML} to multi-level hyperparameters for the flow past RAE2822 airfoil. Top Row: Lift, Bottom Row: Drag. We plot prediction error (Y-axis) vs. Computational Cost (X-axis)}
\label{fig:fpa3}
\end{figure}

Next, we study the sensitivity of prediction errors and in particular, of the gains of the multi-level algorithm over the single-level algorithm, with respect to the three multi-level hyperparameters i.e, model complexity $c_{ml}$, number of training samples $N_0$ (at the coarsest level) and $N_L$ (at the finest level). This sensitivity, for both the lift and the drag, is plotted in figure \ref{fig:fpa3}. We have the following observations from this figure,
\begin{itemize}
    \item Sensitivity to $N_0$: We observe from figure \ref{fig:fpa3} (Left column) that the best gains for the multi-level model arise when a larger number of samples is used at the coarsest level of resolution.
    \item Sensitivity to $N_L$: We observe from figure \ref{fig:fpa3} (Middle column) that the best gains for the multi-level model arise when a few samples are used at the finest level. This is not surprising as increasing the number of samples at the finest level increases the overall cost dramatically, without possibly reducing the prediction error to the same extent.
    \item Sensitivity to $c_{ml}$: We observe from figure \ref{fig:fpa3} (Right column) that the models of intermediate complexity, particularly with $c_{ml} = 1$, provide the best gains over the single-level model as they ensure a balance between accuracy and computational cost, when compared to the inaccurate models of low complexity and costly models of high complexity. 
\end{itemize}
\subsubsection{Uncertainty quantification}
Finally, we consider uncertainty quantification (forward UQ) by computing approximations of the push-forward measure (probability distribution) \eqref{eq:pf1}, with respect to each observable. To this end, we use the generalized version of the ML2MC algorithm (by replacing the multi-level algorithm \ref{alg:ML} with the extended multi-level algorithm \ref{alg:ext_ML}). For comparison, a reference Quasi-Monte Carlo solution is computed with $2000$ Sobol points. 

We compute approximations to push-forward measure with multi-level configurations, shown in Table \ref{tab:MLML_sum} for the lift and the drag. For the sake of comparison, we also compute the measure with the standard single-level Quasi-Monte Carlo algorithm (QMC), the multi-level Quasi-Monte Carlo (MLQMC) algorithm \cite{LyePhD} and the single-level machine learning algorithm (SL2MC) trained with 4, 7, 16, 32, 64, 128, 144, 228, 256 samples for the lift and 4,~8,~16,~32,~110,~128 samples for the drag. For the MC method, we choose the number of samples such that the computational costs of the algorithms are comparable to the cost of the multi-level algorithm.

\begin{figure}[htbp]
    \begin{subfigure}{.45\textwidth}
        \centering
        \includegraphics[width=1\linewidth]{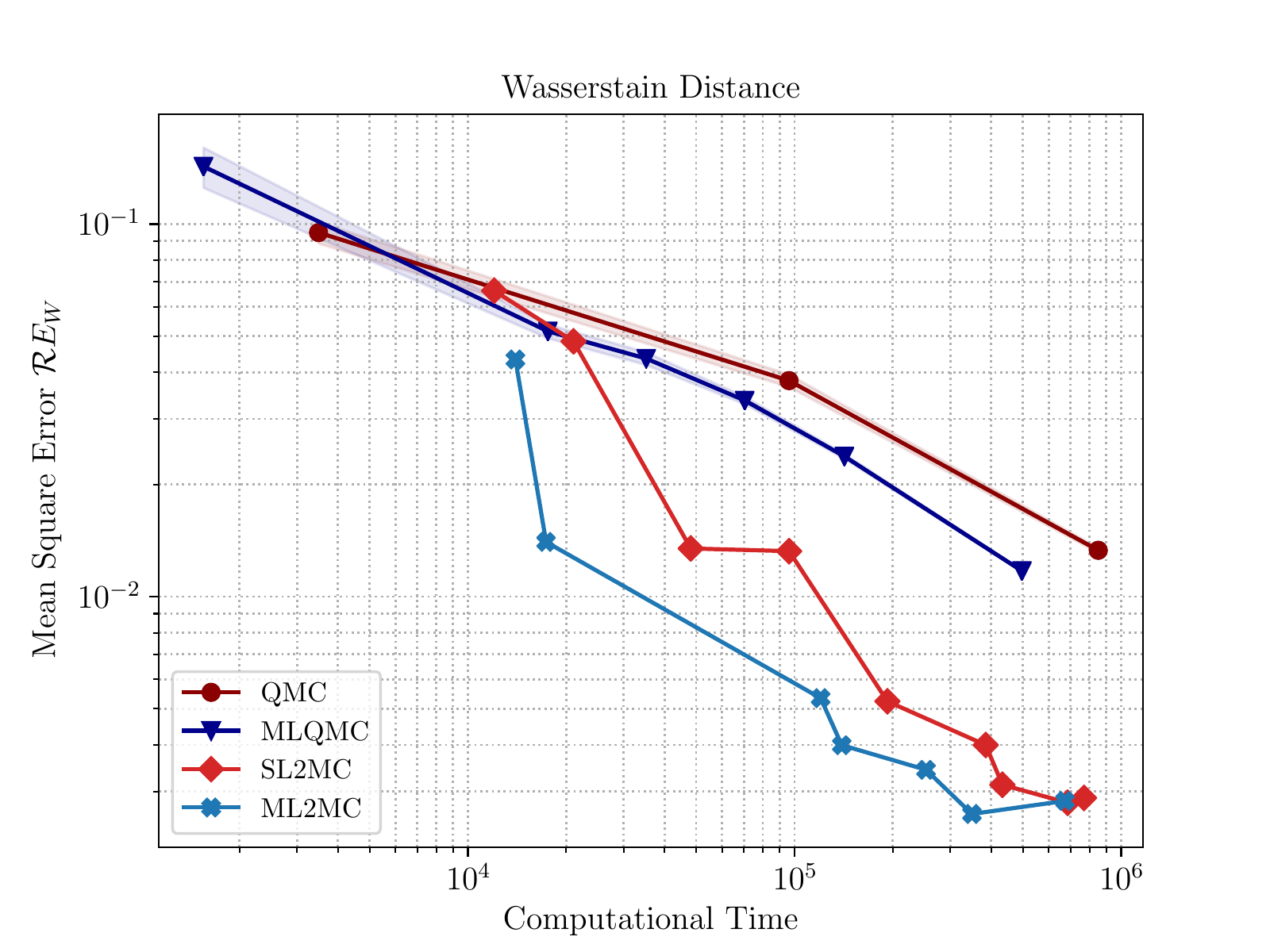}
        \caption{Lift}
    \end{subfigure}
    \begin{subfigure}{.45\textwidth}
        \centering\
        \includegraphics[width=1\linewidth]{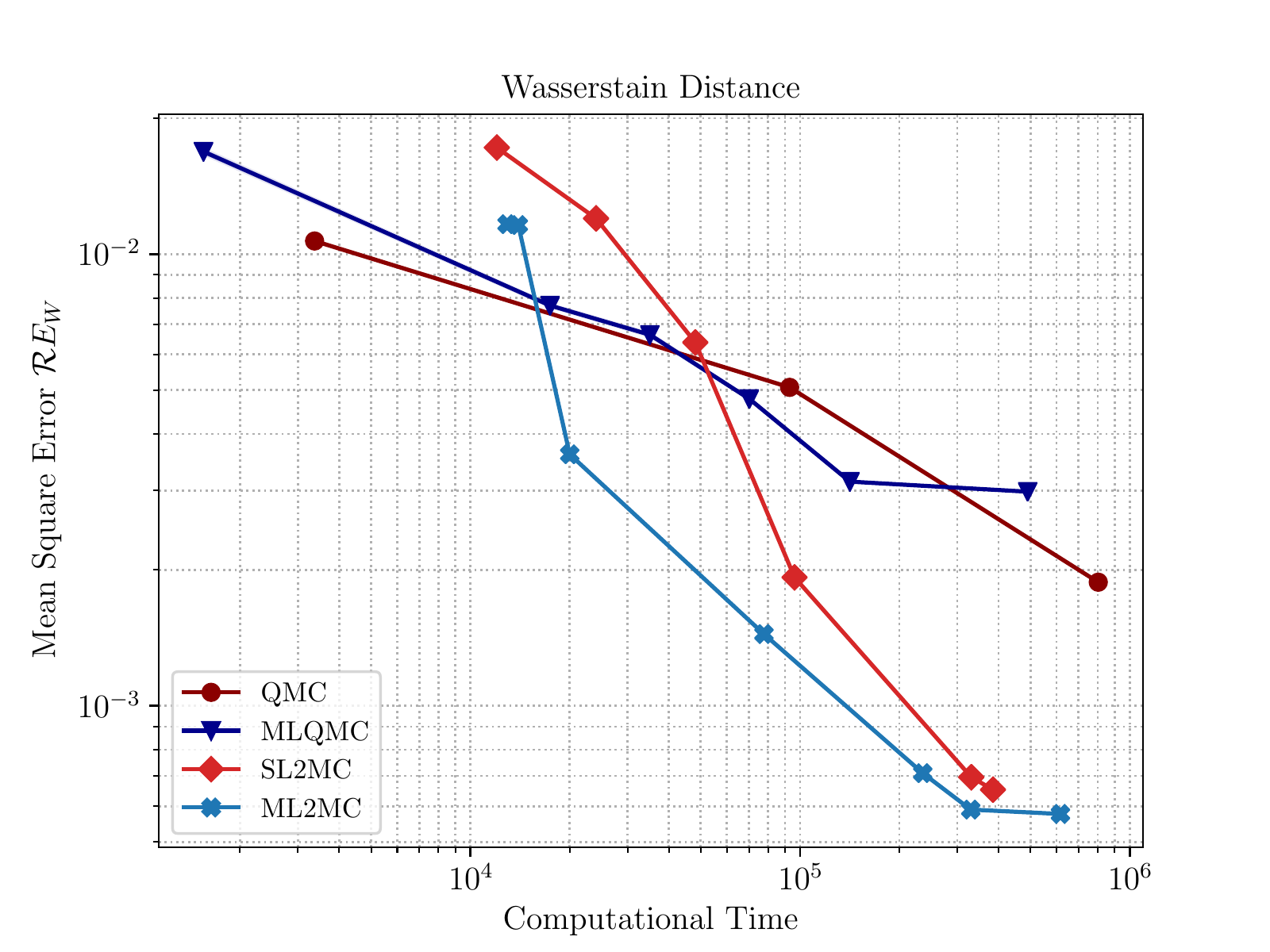}
        \caption{Drag}
    \end{subfigure}
    \caption{Uncertainty quantification for the RAE2822 airfoil. The Wasserstein distance between the approximate push-forward measure \eqref{eq:pf1} computed with standard QMC, multi-level QMC (MLQMC), single-level machine learning (SL2MC) and multi-level machine learning (ML2MC) and a reference push-forward measure for the Lift and the Drag vs computational cost}
    \label{fig:fpa4}
    \end{figure}

The results are shown in figure \ref{fig:fpa4}, where we plot the Wasserstein distance between the approximations and the reference measure, as a function of the computational cost. As seen from the figure, the UQ algorithm, based on the multi-level machine learning algorithm \ref{alg:ext_ML}, clearly and substantially outperforms the standard QMC, multi-level QMC, and the single level machine learning algorithm for both the lift and the drag by providing considerably lower error for the same computational cost. In particular, we obtain an average speedup of $9.2$ (over QMC), $5.4$ (over multi-level QMC) and $1.7$ (over single-level machine learning) for the lift and an average speedup of  $4.8$ (over QMC), $3.5$ (over multi-level QMC) and $1.8$ (over single-level machine learning) for the drag. The maximum speed-ups are even more impressive as they range from $11.4$ (over QMC), $6.4$ (over multi-level QMC) and $4.0$ (over single-level machine learning) for the lift to an average speedup of  $6.8$ (over QMC), $5.2$ (over multi-level QMC) and $3.8$ (over single-level machine learning) for the drag. We would like to remark that the speed-ups can be even more substantial for the moments with respect to the push-forward measure. For instance for the standard deviation, we obtain speedup, on an average over the range of computational costs, of $5.9$ (over QMC), $4.1$ (over multi-level QMC) and $4.1$ (over single-level machine learning) for the lift and an average speedup of  $5.5$ (over QMC), $3.3$ (over multi-level QMC) and $3.0$ (over single-level machine learning SL2MC) for the drag. 
\begin{table}[!t] 
    \centering 
    \renewcommand{\arraystretch}{1.28}
    \footnotesize
    \begin{tabular}{c |c |c |  c |  c| c| c | c | c  } 
        \toprule
        \multicolumn{9}{c}{\bfseries ML2MC Configurations}\\
        \midrule
        \multirow{ 3}{*}{{\bfseries Lift}}
        &\bfseries Samples $N_0$  & 512 & 256 &   2048 &2048& 2048 &2048&2048 \\
        &\bfseries Samples $N_L$ & 4 & 4 &8&32 &64& 92 &92 \\
        &\bfseries Complexity $c_{ml}$ & 0.25  & 1  &  4& 1 &1 &1 &2.31 \\
        \midrule
        \multirow{ 3}{*}{{\bfseries Drag}}
        &\bfseries Samples $N_0$  & 256 & 512 & 512 & 512 & 1024& 1024 & 1024   \\
        &\bfseries Samples $N_L$ & 4 & 4 & 4 & 8 & 64& 92 & 92  \\
        &\bfseries Complexity $c_{ml}$ & 0.25 & 0.25 & 1 & 4 & 1 & 1 & 2.31 \\
        \bottomrule 
    \end{tabular}
    \caption{Different configurations for multilevel machine learning model for approximating the push forward measure wrt Lift and Drag.}
    \label{tab:MLML_sum}
\end{table}

\subsection{Code.} The assembling of the multilevel models and the ensemble training for the selection of the model hyperparameters are performed with a collection of Python scripts, with the support of Keras, Tensorflow and Scikit-learn. The scripts for the generation of the data set for the first numerical experiments and the construction of the multilevel model for both experiments can be downloaded from \url{https://github.com/mroberto166/MultilevelMachineLearning}.

\section{Discussion}
\label{sec:7}
Machine learning, particularly deep learning, algorithms are increasingly popular in the context of scientific computing. One very promising area of application of these algorithms is in the computation of \emph{observables}, corresponding to systems modeled by PDEs. The computation of these observables is very expensive as PDEs have to be simulated for every query i.e, every call of the underlying parameters to observable map \eqref{eq:ptoob}. Instead, following \cite{LMR1}, one can train deep neural networks to provide a surrogate for this parameters to observable map. Although it works well in practice, it was already observed in \cite{LMR1} that finding and training a neural network to approximate the parameters to observable map is very challenging as one can only expect to compute a few training data points (samples), given that the evaluation of each sample involves a very expensive PDE solve. 

We tackle this issue in this paper and present a novel \emph{multi-level} algorithm to significantly increase the accuracy of deep learning algorithms, particularly in the poor data regime. The key idea behind our algorithm is based on the upper bound \eqref{eq:genE2} on the generalization error of a deep learning algorithm for regressing maps. We see from \eqref{eq:genE2} that a significant component of this error is the variance (standard deviation) of the underlying map. Hence, variance reduction techniques can help reduce the generalization error. 

Multi-level methods are examples for a class of variance reduction techniques, termed as control variate methods and are heavily used in the context of uncertainty quantification. We adapt the multi-level idea to machine learning algorithms. The main principle here is to simulate the training data on several mesh resolutions. A large number of cheap (computationally) training samples are used at coarse mesh resolutions whereas only a few computational expensive training samples at generated at fine mesh resolutions, to learn the details i.e differences between successive mesh resolutions. 

We provide theoretical arguments in the form of lemmas \ref{lem:1} and \ref{lem:2} to support our contention that under some reasonable hypothesis on the training process, the multi-level machine-learning algorithm will lead to a decrease in generalization error, when compared to a single-level deep learning algorithm, at the same computational cost. 

In fact, the same design principle works for a very general multi-level machine learning algorithm \ref{alg:ext_ML} that learns the parameters to observable map by a judicious combination of multi-level training on either random points or low-discrepancy sequences, deep neural networks and Gaussian process regressions. 

We test the proposed algorithms on two representative problems. The first is a toy problem of ODEs modeling projectile motion but with a high-dimensional parameter space and highly non-linear parameters to observable map. The second problem is a benchmark uncertain compressible flow past an airfoil. For both problems, we observe that the multi-level algorithm significantly outperforms the single-level machine learning algorithm, resulting in computational gains from half an order to an order of magnitude, in the data poor regime. Moreover, we provide a recipe for finding the set of (multi-level) hyperparameters that result in the highest gains.  

The multi-level algorithm is the basis for a machine-learning algorithm \ref{alg:ML2MC} for forward UQ or uncertainty propagation. Again we observe from the numerical experiments that the multi-level algorithm outperforms all competing algorithms that we tested and provided a computational gain of half an order to an order of magnitude over the standard MC (QMC), MLMC (MLQMC) and single-level machine learning models (such as the deep learning DLMC and DLQMC UQ algorithms of \cite{LMR1}). 

Based on both theoretical and empirical results, we conclude that the multi-level algorithm provides a simple, straightforward to implement and efficient method for improving machine learning algorithms in the context of scientific computing. 

Although we consider only the example of non-linear hyperbolic PDEs, the theory and the algorithms are readily extended to other PDEs such as elliptic and parabolic PDEs. 

The results of this paper can be extended in many different directions, for instance the so-called \emph{multi-fidelity} algorithms \cite{MLFID} can be readily adapted to machine learning and a multi-fidelity algorithm can be designed.

The multi-level algorithm can be extended to learn the whole solution field of the PDE \eqref{eq:ppde}, instead of to just the observable. 

Moreover, one can envisage applying multi-level techniques for problems beyond PDEs, in fact beyond traditional scientific computing. In fact, one could use the multi-level techniques in more traditional applications of machine learning such as image and speech processing, where the multiple levels correspond to different resolutions of the image or the sound file. An interesting application of a related multi-level method in the context of computational chemistry is provided in the recent paper \cite{Har}. Finally, the multi-level machine learning algorithm \ref{alg:ext_ML} will be used in the context of shape optimization and Bayesian inverse problems in future work.
\section*{Acknowledgements.} The research of SM was partially supported by European Research Council Consolidator grant ERCCoG 770880: COMANFLO.


\begin{thebibliography}{99}

\bibitem{ARORA}
Sanjeev Arora, Rong Ge, Behnam Neyshabur, and Yi Zhang. 
\newblock Stronger generalization bounds for
deep nets via a compression approach. 
\newblock {\em In Proceedings of the 35th International Conference on
Machine Learning,} volume 80, pages 254-263. PMLR, Jul 2018.

\bibitem{CAF1}
R. E. Caflisch.
\newblock Monte Carlo and quasi-Monte Carlo methods.
\newblock {\em Acta. Numer.}, 1, 1988, 1-49.

\bibitem{CS1}
F. Cucker and S. Smale.
\newblock On the mathematical foundations of learning.
\newblock {\em Bull. Amer. Math. Soc.,} 39 (1), 2001, 1-49.

\bibitem{Cy1}
G. Cybenko.
\newblock Approximations by superpositions of sigmoidal functions.
\newblock {\em Approximation theory and its applications.,} 9 (3), 1989, 17-28


\bibitem{E1}
W. E and B. Yu.
\newblock The deep Ritz method: a deep learning-based numerical algorithm for solving variational problems.
\newblock {\em Commun. Math. Stat.} 6 (1), 2018, 1-12.

\bibitem{E2}
W. E, J. Han and A. Jentzen.
\newblock Deep learning-based numerical methods for high-dimensional parabolic partial differential equations and backward stochastic differential
equations.
\newblock  {\em Commun. Math. Stat.} 5 (4), 2017, 349-380.

\bibitem{EMW1}
W. E, C. Ma and L. Wu.
\newblock A priori estimates for the generalization error for two-layer neural networks.
 {\em ArXIV preprint}, available from  arXiv:1810.06397, 2018. 




\bibitem{AFOL}
R. Evans et. al.
\newblock De novo structure prediction with deep learning based scoring.
\newblock {\em Google DeepMind working paper}, 2019.


\bibitem{GIL1} M. B. Giles.
\newblock \emph{Multilevel Monte Carlo path simulation.}
\newblock Oper. Res., {\bf 56}:607-617, 2008.

\bibitem{GIL2}
M. B. Giles.
\newblock Multilevel Monte Carlo methods.
\newblock {\em Acta Numer.,} 24 (2015), 259-328.

\bibitem{DLbook}
\newblock I. Goodfellow, Y. Bengio, A. Courville.
\newblock Deep learning.
\newblock{\em MIT press, 2016}

\bibitem{E3}
J. Han. A. Jentzen and W. E.
\newblock Solving high-dimensional partial differential equations using deep learning.
\newblock {\em PNAS}, 115 (34), 2018, 8505-8510.

\bibitem{UMRIDA} C. Hirsch, D. Wunsch, J. Szumbarksi, L. Laniewski-Wollk and J. pons-Prats (editors).
\newblock \emph{Uncertainty management for robust industrial design in aeronautics}
\newblock Notes on numerical fluid mechanics and multidisciplinary design (140), Springer, 2018.

\bibitem{Hein1} S. Heinrich.
\newblock \emph{Multilevel Monte Carlo methods.}
\newblock Large-scale scientific computing, Third international
conference LSSC 2001, Sozopol, Bulgaria, 2001, 
Lecture Notes in Computer Science, Vol {\bf 2170}, 
Springer Verlag (2001), pp. 58-67.



\bibitem{Kor1}
K. Hornik, M. Stinchcombe, and H. White. 
\newblock Multilayer feedforward networks are universal
approximators.
\newblock \emph{ Neural networks,} 2(5), 359-366, 1989.

\bibitem{ADAM} Diederik P. Kingma and Jimmy Lei Ba. 
\newblock Adam: a Method for Stochastic Optimization. 
\newblock {\em International
Conference on Learning Representations}, 1-13, 2015.

\bibitem{Lag1}
I. E. Lagaris, A. Likas and D. I. Fotiadis.
\newblock Artificial neural networks for solving ordinary and partial differential equations.
\newblock {\em IEEE Transactions on Neural Networks}, 9 (5), 1998, 987-1000.

\bibitem{DL-nat}
Y. LeCun, Y. Bengio and G. Hinton.
\newblock Deep learning.
\newblock {\em Nature}, 521, 2015, 436-444.

\bibitem{design_exper}
\newblock  Sacks et al.
\newblock Design and analysis of computer experiments.
\newblock {\em Statistical science (1989)}: 409-423.

\bibitem{LMR1}
\newblock K. O. Lye, S. Mishra, D. Ray. 
\newblock Deep learning observables in computational fluid dynamics. 
\newblock{\em Journal of Computational Physics}:109339 (2020).

\bibitem{LyePhD}
\newblock K. O. Lye.
\newblock Statistical solutions of hyperbolic systems of conservation laws.
\newblock {\em PhD Thesis,} ETH Zurich.

\bibitem{DL_SM1}
S. Mishra.  
\newblock A machine learning framework for data driven acceleration of computations of differential equations, 
\newblock {\it Math. in Engg.}, 1 (1), 2018, 118-146. 

\bibitem{MR1}
S. Mishra and K. Rusch.  
\newblock Enhancing accuracy of deep learning algorithms by training with low-discrepancy sequences.
\newblock {\it Preprint,}, 2020. Available from arXiv:2005.12564.

\bibitem{MS1}
S. Mishra and C. Schwab.
\newblock Sparse tensor multi-level Monte Carlo finite volume methods for hyperbolic conservation laws with random initial data.
\newblock {\it Math. Comput.}, 81(180), 2012, 1979--2018.

\bibitem{MSS1}
S. Mishra, Ch. Schwab and J. {\v S}ukys.
\newblock Multi-level Monte Carlo finite volume methods for nonlinear systems of conservation laws in multi-dimensions.
\newblock \textit{J. Comput. Phys} 231 (8), 2012, 3365--3388.

\bibitem{MSDR}
T. De Ryck, S. Mishra, R. Deep.
\newblock On the approximation of rough functions with deep neural networks.
\newblock \textit{Preprint}, 2020. Available from arXiv:1912.06732.


\bibitem{JR1}
T. P. Miyanawala and R. K, Jaiman.
\newblock An efficient deep learning technique for the Navier-Stokes equations: application to unsteady wake flow dynamics.
\newblock {Preprint}, 2017, available from arXiv :1710.09099v2.

\bibitem{NEYS1}
Behnam Neyshabur, Zhiyuan Li, Srinadh Bhojanapalli, Yann LeCun, and Nathan Srebro. 
\newblock Towards understanding the role of over-parametrization in generalization of neural networks. 
{\em arXiv
preprint} arXiv:1805.12076, 2018.

\bibitem{MLFID}
B. Peherstorfer, K. Willcox, and M. Gunzburger.
\newblock Survey of multifidelity methods in uncertainty propagation, inference, and optimization. 
{\em SIAM Rev.} 60 (2018), no. 3, 550–591.

\bibitem{Kar1}
M. Raissi and G. E. Karniadakis.
\newblock Hidden physics models: machine learning of nonlinear partial differential equations.
\newblock  {\em J. Comput. Phys.,} 357, 2018, 125-141.

\bibitem{Kar2}
M. Raissi, A. Yazdani and G. E. Karniadakis.
\newblock Hidden fluid mechanics: Learning velocity and pressure fields from flow visualizations.
\newblock  {Science}, 367.6481 (2020), 1026-1030

\bibitem{Kar3}
L. Lu, P. Jin and G. E. Karniadakis.
\newblock DeepONet: Learning nonlinear operators for identifying differential equations based on the universal approximation theorem of operators
\newblock  {Preprint}, 2019. Available from arXiv:1910.03193


\bibitem{GP}
\newblock C.E Rasmussen. 
\newblock Gaussian processes in machine learning. 
\newblock {\em Summer School on Machine Learning}, Springer, Berlin, Heidelberg, 2003.

\bibitem{DR1}
D. Ray and J. S, Hesthaven.
\newblock An artificial neural network as a troubled cell indicator.
\newblock {\em J. Comput. Phys.,}  367, 2018, 166-191

\bibitem{RCFM}
 D. Ray, P. Chandrasekhar, U. S. Fjordholm and S. Mishra.
 \newblock Entropy stable scheme on two-dimensional unstructured grids for Euler equations, 
 \newblock {\it Commun. Comput. Phys.,} 19 (5), 2016, 1111-1140.
 
 
 
 \bibitem{SG}
S. Ruder.
\newblock An overview of gradient descent optimization algorithms.
\newblock {Preprint}, 2017, available from arXiv.1609.04747v2.

\bibitem{MLbook}
Shai Shalev-Shwartz and Shai Ben-David.
\newblock{\em Understanding machine learning: From theory to
algorithms.} 
\newblock Cambridge university press, 2014.

\bibitem{INC}
J. Tompson, K. Schlachter, P. Sprechmann and K. Perlin.
\newblock Accelarating Eulerian fluid simulation with convolutional networks.
\newblock {\em Preprint}, 2017. Available from arXiv:1607.03597v6.

\bibitem{ROMbook}
A. Quateroni, A. Manzoni and F. Negri.
\newblock \emph{Reduced basis methods for partial differential equations: an introduction},
\newblock Springer Verlag (2015).

\bibitem{YAR1}
D. Yarotsky. 
\newblock Error bounds for approximations with deep ReLU networks. 
\newblock {\em Neural Networks}, 
94, 2017, 103-114

\bibitem{Har}
P. Zaspel, B. Huang, H. Harbrecht and O. Anotole von Lillenfeld.
\newblock Boosting quantum machine learning with multi-level combination technique: Pople diagrams revisited.
\newblock {Preprint,} available as Arxiv1808.02799v2.

\bibitem{py_wass}
Statistical functions (scipy.stats).
\newblock Python Library.
\newblock \url{https://docs.scipy.org/doc/scipy/reference/stats.html}

\end{thebibliography}
\end{document}